\scrollmode

\documentclass[a4paper,11pt]{article}
\usepackage{amssymb,amsmath}
\usepackage{hyperref}

\def\YM{Y\!M}
\def\R{{\Bbb R}}

\def\N{{\Bbb N}}
\def\C{{\Bbb C}}
\def\i{\mathbf{i}}
\def\eas{\begin{eqnarray*}}
\def\eeas{\end{eqnarray*}}
\def\lll{\lefteqn}
\def\nn{\nonumber}
\def\eq#1{\begin{equation}\label{#1}}
\def\eeq{\end{equation}}
\def\ea#1{\begin{eqnarray}\label{#1}}
\def\eea{\end{eqnarray}}
\def\la{\label}
\def\re#1{(\ref{#1})}
\def\dist{\mathop{\mbox{\rm dist}}}

\def\a{{\alpha}}

\def\d{\partial}

\def\p{{\varphi}}
\def\dx{\,dx}

\def\del{{\delta}}
\def\eps{{\varepsilon}}

\def\qed{\hfill$\Box$\\\strut}

\def\loc{{\mbox{\rm \tiny loc}}}

\def\st{\,{*}\,}
\def\<{\langle}
\def\>{\rangle}
\def\Hat{\widehat}
\def\Kreuz{\mathop{\mbox{\Large\#}}}
\newcommand{\wto}{\rightharpoondown}
\newcommand{\wsto}{\overset{\raisebox{-1ex}{\scriptsize
      $*$}}{\rightharpoondown}}
\newcommand{\HdR}{H_{\rm dR}}

\newtheorem{lemma}{Lemma}[section]

\newtheorem{corollary}[lemma]{Corollary}
\newtheorem{proposition}[lemma]{Proposition}
\newtheorem{theorem}[lemma]{Theorem}
\newtheorem{remark}[lemma]{Remark}
\def\bib#1{\bibitem[#1]{#1}}
\def\gg{{\mathfrak g}}
\def\suk{{\mathfrak{su}}(k)}
\def\ddt{\frac{d}{dt}_{|t=0}}
\def\da#1{d_A^{\wedge{#1}}}
\def\das#1{d_A^{*\wedge{#1}}}
\def\dast#1{d_{A_t}^{*\wedge{#1}}}
\def\dasi#1{d_{A_i}^{*\wedge{#1}}}

\def\dn#1{d^{\wedge{#1}}}

\def\dnx#1{d^{(*)\wedge{#1}}}
\def\tf#1#2{{\textstyle\frac{#1}{#2}}}
\def\Om{{\Omega}}

\def\edge{\hspace{0.15em}\mbox{\LARGE$\llcorner$}\hspace{0.05em}}

\def\XXint#1#2#3{{\setbox0=\hbox{$#1{#2#3}{\int}$}
                 \vcenter{\hbox{$#2#3$}}\kern-.5\wd0}}
\def\Xint#1{\mathchoice
              {\XXint\displaystyle\textstyle{#1}}
              {\XXint\textstyle\scriptstyle{#1}}
              {\XXint\scriptstyle\scriptscriptstyle{#1}}
              {\XXint\scriptscriptstyle\scriptscriptstyle{#1}}
              \!\int}
\def\mint{\Xint-}

\begin{document}
\title{Minimizers of\\ higher order gauge invariant functionals}
\author{Andreas Gastel \and Christoph Scheven}
\maketitle

\numberwithin{equation}{section}
\setcounter{section}{0}
\setlength{\parindent}{0em}

\begin{abstract}\noindent
  We introduce higher order variants of the Yang-Mills functional
  that involve $(n-2)$th order derivatives of the curvature. We prove
  coercivity and smoothness of critical points in Uhlenbeck
  gauge in dimensions $\mathrm{dim}M\le 2n$. 
  These results are then used to establish the existence of
  smooth minimizers on a given principal bundle $P\to M$ for subcritical
  dimensions $\mathrm{dim}M<2n$. In the case of critical dimension
  $\mathrm{dim}M=2n$ we construct a minimizer on a bundle which might
  differ from the prescribed one, but has the same Chern classes 
  $c_1,\ldots,c_{n-1}$. A key result is a removable singularity theorem for
  bundles carrying a $W^{n-1,2}$-connection. This generalizes a recent
  result by Petrache and Rivi\`ere.\\\strut


\end{abstract}

\section{Introduction}

One of the limitations of elliptic Yang-Mills theory as a variational problem
is that many of its features work only if the basis of the bundle has dimension
$\le4$. Given a principal $G$-bundle $P\to M$ over a compact manifold $M$ with
compact structure group $G$, we consider $G$-connections $D_0+A$ on it.
(We will explain basic concepts of gauge theory in 
Section~\ref{sec:gauge-theory}.)
The squared curvature integral
\[
  \YM(A):=\frac12\int_M|F_A|^2\,dx
\]
is the {\em Yang-Mills functional\/}. Its variational theory is based on
fundamental contributions by Uhlenbeck \cite{Uh1} \cite{Uh2} and Sedlacek
\cite{Se} from the 1980s. If $\dim M\le3$, there is a smooth minimizer
of $\YM$ on any given bundle. If $\dim M=4$, which is the critical and most
natural case, then things are getting more interesting. A minimizing sequence
$(A_j)_{j\in\N}$ for $\YM$ in a given bundle has to be looked at in a good 
gauge (roughly a good choice of bundle coordinates) in order to control 
Sobolev norms of $A_j$ well enough to find a subsequence converging weakly
to some sort of minimizer $A$. The latter turns out to be smooth on most
of $M$, but may initially have point singularities. Seeing those as possible
singularities of the bundle itself, they can be removed to get a smooth
minimizer $A$, but possibly on a different bundle, which is redefined even
topologically in the singular points. The bundle in which to find the minimizer,
however, is not arbitrary, since some of its topological invariants coincide
with those of the bundle we started with.

It is the gauge choice for the minimizing sequence that may fail in dimensions
$\dim M\ge5$. This can be fixed partially, since the gauge theorem from
\cite{Uh2} works just as well on $2n$-dimensional manifolds for
\[
  \YM^n(A):=\frac1n\int_M|F_A|^n\,dx
\]
instead of $\YM$. Therefore, there may be hope that the variational approach 
we described for $\YM$ on $4$-dimensional manifolds can be modified to
work for $\YM^n$ on $2n$-dimensional manifolds. Part of this program has
been worked out, namely \cite{Uh2} works for the gauge choice, and
Isobe \cite{Is} has worked out some local regularity theorem which can be
applied to minimizers of $\YM^n$. However, we think that there are two points
where minimizing $\YM^n$ in $2n$ dimensions is not as good as minimizing
$\YM$ in $4$ dimensions.
\begin{enumerate}\renewcommand{\labelenumi}{(\arabic{enumi})}
\item Minimizers cannot be expected to be smooth, since the functional is
degenerate. The best we can hope for, and that is essentially what has been
proven in \cite{Is}, is $C^{1,\alpha}$-regularity (except maybe for
point singularities).
\item Uhlenbeck's singularity removal theorem from \cite{Uh1} has been improved
significantly by Petrache and Rivi\`ere \cite{PR} by removing the assumption
that $A$ is Yang-Mills. In fact, the existence of {\em any\/} 
$W^{1,2}$-connection on a bundle over $M$ minus one point implies that the 
bundle can be continued to give a bundle over all of $M$, provided that
$\dim M\le4$. It is our impression that the arguments from \cite{PR}
do not carry over to $W^{1,n}$-connections for $\dim M=2n$.
\end{enumerate}

The starting point for our paper is the following. {\em Both problems (1)
and (2) do not occur if we work in $W^{n-1,2}$ instead of $W^{1,n}$ for
$\dim M\le 2n$.\/} As we shall prove in Theorem \ref{RemSing}, the existence
of a $W^{n-1,2}$-connection on a $2n$-dimensional bundle with one fibre
missing implies that the singularity can be removed from the bundle. This
directly generalizes the corresponding theorem from \cite{PR} and therefore
helps to get around problem (2). To handle the issue mentioned in (1), we
have to work with {\em nondegenerate\/} functionals that control the
$W^{n-1,2}$-norm (and hence also the $W^{1,n}$-norm). Therefore 
let us try to write 
down functionals that do the job and can then be used to replace $\YM^n$.

A word on gauge invariance first, which is one of the basic features that
make the Yang-Mills functional interesting.
For any \emph{gauge transformation}, i.e. any sufficiently regular
equivariant map
$u:P\to G$, we have $\YM(A)=\YM(u^*A)$, because $F_A$ transforms like
\[
  F_{u^*A}=u^{-1}F_Au,
\]
and $u$ as well as $u^{-1}$ act by isometries. One of the issues we have to
deal with is finding {\em higher order\/} functionals that show the same
gauge invariance.

There has been some work in this direction, actually, at least for the $n=3$
case. In \cite{BU}, Bejan and Urakawa have defined the Bi-Yang-Mills functional
\[
  \YM_2(A):=\frac12\int_M|d_A^*F_A|^2\,dx,
\]
where here $d_A^*$ is the covariant exterior co-differential. It is 
gauge-invariant, since the Euler-Lagrange equation $d_A^*F_A=0$ of $\YM$ 
must be gauge-invariant. In \cite{IIU1} and \cite{IIU2},
some more of the basic properties of $\YM_2$ are explored.

For our purpose of constructing minimizers, the functional $\YM_2$ is not quite suitable,
since it does not control any $L^p$-norm of $|F_A|$, which is bad news
for controlling minimizing sequences. But we can add $\YM^3$ to it, arriving
at
\[
  Y_3(A):=\int_M(|d_A^*F_A|^2+|F_A|^3)\,dx,
\]
for which, as will follow from our results, all of the Uhlenbeck-Sedlacek
program described above can be performed similarly in dimensions $\le6$.
Philosophically, $Y_3(A)$ can control the $W^{1,2}$-norm of $F_A$ since it
obviously controls $d_A^*F_A$ in $L^2$, and $d_AF_A$ is always $0$ by
Bianchi's identity. Hodge theory says that $DF_A$ is controlled once you
can control $dF_A$ and $d^*F_A$, so what we need is some ``nonlinear
variant'' of Hodge theory. We cannot work with $|d^*F_A|$ 
directly, because it is not gauge invariant.

We go on constructing higher order gauge invariant functionals inductively.
First derivatives $d_A^*F_A$ should be controlled by norms of $d_Ad_A^*F_A$
and $d_A^*d_A^*F_A$, both of which have gauge invariant norms, see
Section~\ref{sec:diff-forms}. But $d_A^*d_A^*F_A=-*[F_A,*F_A]$ happens to
be of lower order, so $d_Ad_A^*F_A$ alone should be enough to control
first derivatives $d_A^*F_A$, and hence second derivatives of $F_A$. The
functional
\[
  Y_4(A):=\int_M(|d_Ad_A^*F_A|^2+|F_A|^4)\,dx,
\]
does the job, for $\dim M\le8$. Of course, we can iterate our 
considerations and find functionals suitable for our program. Abbreviating
\[
  \das{n}:=\left\{\begin{array}{ll}
    (d_Ad_A^*)^{n/2} & \mbox{ if $n$ even,}\\
    d_A^*(d_Ad_A^*)^{(n-1)/2} & \mbox{ if $n$ is odd,}
  \end{array}\right.
\]
we define
\[
  Y_n(A):=\int_M(|\das{n-2}F_A|^2+|F_A|^n)\,dx.
\]
These are gauge invariant, and scaling invariant if $\dim M=2n$. Moreover,
they will turn out to be coercive when put in the gauge found by Uhlenbeck
\cite{Uh2}. And, of course, being quadratic in the highest order, they are 
then nondegenerate, which opens the possibility of proving $C^\infty$ for
minimizers instead of $C^{1,\a}$, thus addressing problem (1). It looks like
we have found good candidates for functionals to look at.

From another point of view, the functionals $Y_n$ may be not the best choice.
They are no perturbations of the original Yang-Mills functional. We may wish to
minimize ``Yang-Mills plus something of higher order'' and even think of
that ``something'' being multiplied by some small $\eps>0$. The higher
order terms we described so far do need the $|F_A|^n$-term in order to be
coercive. But if we are prepared to leave the realm of exterior forms, we
can proceed. Instead of using exterior derivatives $d_A$ and $d_A^*$, we can
try to use other combinations of exterior partial derivatives. It turns out
that the norm of the {\em full\/} covariant derivative $|D_AF_A|$ is also
gauge invariant, and so are its iterates $|D_A^kF_A|$. Using these, we come up
with a second sequence of functionals, this time ``perturbations of $\YM$'',
which read
\[
  Z_n(A):=\int_M(|D_A^{n-2}F_A|^2+|F_A|^2)\,dx.
\]
Both sorts of functionals have their advantages, and it turns out that
they can be estimated against each other and against the (squared) ``nonlinear
$W^{n-2,2}$-norms'' of $F_A$, which are built like the usual Sobolev norms, but
using $D_A$ instead of $D$.

For the functionals $Y_n$ and $Z_n$, we will prove the following results 
concerning the existence and regularity of minimizers.
\begin{description}
\item[Section~\ref{sec:WA}.] 
\emph{Global ``coercivity''} in the sense that $Y_n$ and $Z_n$ either
control the nonlinear $W^{n-2,2}$-norms of $F_A$ mentioned above. This
is not real coercivity since it works only with gauge invariant quantities and
controls $F_A$ instead of $A$ itself. 
The exact statement is Theorem~\ref{WA-estimates}.

\item[Section~\ref{sec:uhlenbeck}.] 
\emph{Local coercivity} in Uhlenbeck gauge in $\le 2n$ dimensions. 
This is essential for extracting weak limits from minimizing sequences, and is 
performed in Theorem~\ref{uhl_n}.

\item[Section~\ref{sec:rem}.] 
\emph{Removability of point singularities of bundles. }As remarked above,
a point singularity of a bundle can be removed once we know the existence of
a $W^{n-1,2}$-connection around the missing point. Even better, we need
only a connection for which $Y_n$ or $Z_n$ is finite. This key result helps us
removing point singularities from the minimizing connection constructed
below, which by construction is in $W^{n-1,2}$. See Theorem~\ref{RemSing}.

\item[Section~\ref{sec:euler}.] 
Weak formulations of the \emph{Euler-Lagrange equations}.

\item[Section~\ref{sec:regularity}.] 
\emph{Smoothness of weak solutions} of the Euler-Lagrange equations, again
for $\dim M\le 2n$. 
This regularity result can be found in 
Theorem~\ref{thm:regularity}.
 
\item[Section~\ref{sec:existence-critical}.] 
\emph{Existence of minimizers in the critical dimension} $\dim M=2n$.
As mentioned before, the choice of Uhlenbeck
gauges -- which is necessary to overcome the lack of coercivity of the
functional -- can be achieved uniformly for the minimizing sequence
only away from finitely many points. This results in finitely many
singularities that might develop in the bundle during the minimizing
process. However, the removable singularities theorem helps us to
remove the singularities of the bundle, and then, using our regularity theorem,
also the singularities of the minimizer. This minimizer, singularities
having been removed, lives on a new
bundle that might
differ from the prescribed bundle. But it still has the same 
Chern classes $c_1(P),\ldots,c_{n-1}(P)$ as the original bundle.
For the detailed statement of the result, we refer to 
Theorem~\ref{thm:existence}. 

\item[Section~\ref{sec:existence-subcritical}.]  
\emph{Existence of minimizers in subcritical dimensions} $\dim M<2n$.
In this case we can start with an arbitrary principal bundle with any
compact structure group $G$ and can construct a minimizer on the given
bundle. Moreover, the constructed minimizing connection is smooth by
the regularity theorem.
The existence result is Theorem~\ref{thm:low_existence}.
\end{description}

\section{Basics}\label{sec:basics}

\subsection{Basic facts on gauge theory}
\label{sec:gauge-theory}

In this section we briefly recall those facts on
connections on principal bundles that will be relevant for the present article. For a more
thorough exposition of the theory, we refer to \cite[App. A]{We}.

Throughout this paper, we fix a smooth compact Riemannian
manifold $M$ of dimension $m:=\dim M\le 2n$ and a compact Lie group
$G$, the Lie algebra of which will be denoted by $\gg$. 
A \emph{principal bundle} $\pi:P\to M$ over $M$ with
structure group $G$ can be described by an open cover
$\{U_\alpha\}_{\alpha=1}^L$ of $M$ and local trivializations
$\phi_\alpha:\pi^{-1}(U_\alpha)\to U_\alpha\times G$. 
The trivializations give rise to transition functions  
$\phi_{\alpha\beta}:U_\alpha\cap
U_\beta\to G$ defined by $\phi_{\alpha\beta}\phi_\beta=\phi_\alpha$ 
for all parameters $1\le \alpha,\beta\le L$ with $U_\alpha\cap
U_\beta\neq\varnothing$. From the definition of the transition
functions, it is immediate to check the \emph{cocycle conditions}
\begin{equation}
  \label{cocyle-conditions}
  \phi_{\alpha\alpha}\equiv1
  \qquad\mbox{and}\qquad
  \phi_{\alpha\beta}\phi_{\beta\gamma}
  =
  \phi_{\alpha\gamma}
  \mbox{\quad on }U_\alpha\cap U_\beta\cap U_\gamma,
\end{equation}
provided $U_\alpha\cap U_\beta\cap
U_\gamma\neq\varnothing$. Conversely, any set of smooth functions
$\{\phi_{\alpha\beta}\}_{\alpha,\beta}$ that satisfies the cocycle
conditions~(\ref{cocyle-conditions}) defines a principal $G$-bundle
with transition functions $\phi_{\alpha\beta}$ relative to
the open cover $\{U_\alpha\}$. 

  A \emph{gauge transformation} on $P$ is an equivariant smooth map
  $u:P\to G$. Using the trivializations $\phi_\alpha$ of the bundle
  $P$, the gauge transformation can alternatively be characterized by
  its localizations $u_\alpha:U_\alpha\to G$, which are related by
  the transition identity
  \begin{equation*}
    u_\beta=\phi^{-1}_{\alpha\beta}u_\alpha\phi_{\alpha\beta}
    \qquad\mbox{on }U_\alpha\cap U_\beta.
  \end{equation*}

  A smooth connection $D_0+A$ on the principal bundle $P$ is formally an
  element of the space $D_0+C^\infty(M,T^*M\otimes\gg_P)$,
  where $D_0$ is a fixed smooth reference connection on $P$ and $\gg_P$
  denotes the associated $\gg$-bundle, cf. \cite[App. B]{We}.  
  See \cite[App. A]{We} for the
  precise definition. For most of the present article however,
  it will be sufficient to think of $A$ as of the entity of its 
  localizations
  $(\phi_\alpha)_*A=A_\alpha\in C^\infty(U_\alpha,T^*M\otimes\gg)$
  subject to the trivializations $\phi_\alpha:\pi^{-1}(U_\alpha)\to U_\alpha\times G$ of the
  principal bundle, which are given by
  \begin{equation*}
    A_\alpha(d_p\pi(v))=\phi_\alpha(p)
    A(d_p\pi(v))\phi_\alpha^{-1}+d\phi_\alpha(v) \phi_\alpha^{-1}(p)
  \end{equation*}
   for all $p\in\pi^{-1}(U_\alpha)$ and $v\in T_pP$.
   In fact, the set $\{A_\alpha\}_{\alpha=1}^L$
  contains the same information as $A$. 
  The localizations $A_\alpha$,  $1\le \alpha\le L$ are related by the identity
  \begin{equation}
    \label{A-transform}
    A_\beta
    =
    \phi_{\alpha\beta}^{-1}A_\alpha\phi_{\alpha\beta}
    +
    \phi_{\alpha\beta}^{-1}d\phi_{\alpha\beta}
    \qquad\mbox{on }U_\alpha\cap U_\beta
  \end{equation}
  for all parameters $\alpha,\beta$ for which the latter set is
  non-empty. On the other hand, if we are given a set of smooth
  local connections $A_{\alpha}\in C^\infty(U_\alpha,T^*M\otimes\gg)$, 
  $\alpha=1,\ldots,L$ that satisfy the compatibility
 conditions~(\ref{A-transform}) with respect to the transition
 functions $\phi_{\alpha\beta}$ of the bundle $P$, 
 we can find a smooth connection $A$
 on $P$ with trivializations $(\phi_\alpha)_*A=A_\alpha$. 

 The gauge transformation $u$ acts on the connection $A$ via
 \begin{equation*}
   u^*A:= u^{-1} A u+u^{-1} du.
 \end{equation*}
 The \emph{curvature} of a connection $A$ is given by
  \begin{equation*}
   F_A:= dA +A\wedge A.
 \end{equation*}
  It is well-known that the curvature is gauge-equivariant in the
  sense that for every gauge transformation $u:P\to G$ there holds
  \begin{equation*}
    F_{u^*A}=u^{-1}F_A u.
  \end{equation*}

 Finally, we note that Sobolev spaces of connections can be defined by
 \begin{equation*}
   \mathcal{A}^{k,p}(P):= D_0+W^{k,p}(M,T^*M\otimes \gg_P),
 \end{equation*}
 where $D_0$ is a fixed smooth reference connection on $P$ and $\gg_P$
 denotes the associated $\gg$-bundle, cf. \cite[App. B]{We}. Locally,
 $W^{n-1,2}$-connections on $U_\alpha$ are represented by
 $A_\alpha\in W^{n-1,2}(U_\alpha,T^*M\otimes\gg)$. 
 Accordingly, we will consider local gauge transformations of class
 $W^{n,2}$, in other words, maps $u_\alpha\in
 W^{n,2}(U_\alpha,G)$.

\subsection{Calculations with differential forms}
\label{sec:diff-forms}

We now consider $\gg$-valued differential forms on a coordinate chart
$U\subset\R^m$. 

For a $\gg$-valued $k$-form $A$ and a $\gg$-valued $\ell$-form $B$, we
introduce the abbreviation 
\begin{equation*}
  [A,B]:=A\wedge B-(-1)^{k\ell} B\wedge A.
\end{equation*}
With this notation, we define in the case $k=1$
\begin{equation*}
  d_A B= dB+[A,B]
  \qquad\mbox{and}\qquad
  d^*_AB:= d^*B+(-1)^{m+1}{*}[A,{*}B].
\end{equation*}

We introduce the following notations for higher order exterior
derivatives of differential forms. For 1-forms $C\in
W^{k,1}(U,\wedge^1\R^m\otimes\gg)$, respectively $2$-forms 
$B\in
W^{k,1}(U,\wedge^2\R^m\otimes\gg)$, where $k\in\N$, we use the notations
\begin{equation*}
  d^{\wedge k}C:=
  \begin{cases}
    \underbrace{d d^*\cdots d}_{k} C,&\mbox{ $k$ odd,}\\
    \underbrace{d^*\cdots d}_{k} C,&\mbox{ $k$ even},
  \end{cases}
  \qquad\quad
  d^{*\wedge k}B:=
  \begin{cases}
    \underbrace{d^* d\cdots d^*}_{k} B,&\mbox{ $k$ odd,}\\
    \underbrace{d\cdots d^*}_{k} B,&\mbox{ $k$ even}.
  \end{cases}
\end{equation*}
Similarly, for a connection $A$ we introduce nonlinear versions 
\begin{equation*}
  d^{\wedge k}_AC:=
  \begin{cases}
    \underbrace{d_A d_A^*\cdots d_A}_{k} C,&\mbox{ $k$ odd,}\\
    \underbrace{d_A^*\cdots d_A}_{k} C,&\mbox{ $k$ even},
  \end{cases}
  \qquad
  d^{*\wedge k}_AB:=
  \begin{cases}
    \underbrace{d_A^* d_A\cdots d_A^*}_{k} B,&\mbox{ $k$ odd,}\\
    \underbrace{d_A\cdots d_A^*}_{k} B,&\mbox{ $k$ even}.
  \end{cases}
\end{equation*}
We apply this in particular to $B=F_A$ to define higher order exterior
derivatives $d^{*\wedge k}_AF_A$ of the curvature.

Furthermore, we define higher order total derivatives $D_A^kF_A$ by
exploiting the fact that the connection $A$ 
induces covariant derivatives on vector bundles associated with the
principal bundle $P$. More precisely, $D_A^kF_A$ is a section of the bundle
$\otimes^kT^*M\otimes\wedge^2T^*M\otimes\gg_P$ and is defined
inductively by
\begin{equation*}
  D_A^kF_A:=D_A(D_A^{k-1}F_A)
  \qquad\mbox{for }k\in\N.
\end{equation*}
Here, the first $D_A$ on the right-hand side
denotes the covariant derivative on
$\otimes^{k-1}T^*M\otimes \wedge^2T^*M\otimes\gg_P$ that is induced
by the Levi-Civita connection on $TM$ and the connection $A$ on
$P$.

By $D_A^*$ we denote the formal adjoint of $D_A$. For later reference,
we remark the existence of constants depending on the bundle, such that
\eq{int0}
  |d_AB|\le C|D_AB|,\qquad
  |d_A^*B|+|D_A^*B|\le C(|D_AB|+|B|)
\eeq
for all forms $B$ in the above bundle.

The above definitions of derivatives of $F_A$ provide us with
two classes of higher order functionals, namely
\begin{equation*}
    Y_n(A):=\int_M(|\das{n-2}F_A|^2+|F_A|^n)\,dx,
\end{equation*}
and
\begin{equation*}
    Z_n(A):=\int_M(|D_A^{n-2}F_A|^2+|F_A|^2)\,dx.
\end{equation*}

For our purposes, it is crucial that both types of functionals are gauge 
invariant. This is a consequence of the following three lemmas. 

\begin{lemma} Let $B_A$ be some $\gg$-valued $2$-form that transforms according
to 
\[
  B_{u^*A}=u^{-1}B_Au.
\]
Then we also have
\[
  d_{u^*A}^*B_{u^*A}=u^{-1}d_A^*B_Au.
\]
\end{lemma}

{\bf Proof.} We compute, using the gauge equivariance of $B_A$ and $0=d(uu^{-1})
=du\,u^{-1}+u\,d(u^{-1})$,

\begin{align*}
 &\st d\st B_{u^*A}=\st d\st(u^{-1}B_Au)\\
  &\quad=\st (u^{-1}d(\st B_A)u+d(u^{-1})\wedge(\st B_A)u
  +(-1)^{m-2}u^{-1}(\st B_A)\wedge du)\\
  &\quad=u^{-1}\st d(\st B_A)u-\st (u^{-1}du\wedge u^{-1}(\st B_A)u
  -(-1)^{m-2}u^{-1}(\st B_A)u\wedge u^{-1}du)\\
  &\quad=u^{-1}\st d\st B_Au-\st [u^{-1}du,\st (u^{-1}B_Au)]
\end{align*}
and    
\begin{align*}
  \st[u^* A,\st B_{u^*A}]
  &=\st [u^{-1}Au+u^{-1}du,\st (u^{-1}B_Au)]\\
  &=u^{-1}\st [A,\st B_A]u+\st [u^{-1}du,\st (u^{-1}B_Au)].
\end{align*}
Thus we have
\[
  d^*_{u^*A}B_{u^*A}=-(-1)^mu^{-1}(\st d\st B_A+\st [A,\st B_A])u=u^{-1}d_A^*B_Au
\]
as claimed.\qed

\begin{lemma}\label{lem:trafo-d}
Let $C_A$ be some $\gg$-valued $1$-form that transforms according
to 
\[
  C_{u^*A}=u^{-1}C_Au.
\]
Then we also have
\[
  d_{u^*A}C_{u^*A}=u^{-1}d_AC_Au.
\]
\end{lemma}

{\bf Proof.} As above,
\begin{align*}
  dC_{u^*A}&=d(u^{-1}C_Au)\\
  &=u^{-1}dC_Au+d(u^{-1})\wedge C_Au-u^{-1}C_A\wedge du\\
  &=u^{-1}dC_Au-u^{-1}du\wedge u^{-1}C_Au-u^{-1}C_Au\wedge u^{-1}du\\
  &=u^{-1}dC_Au-[u^{-1}du,u^{-1}C_Au],
\end{align*}
as well as
\begin{align*}
  {}[u^*A,u^{-1}C_Au]
  &=[u^{-1}Au+u^{-1}du,u^{-1}C_Au]\\
  &=u^{-1}[A,C_A]u+[u^{-1}du,u^{-1}C_Au].
\end{align*}
This implies
\[
  d_{u^*A}C_{u^*A}=dC_{u^*A}+[u^*A,u^{-1}C_Au]=u^{-1}d_AC_Au
\]
as desired.\qed

\begin{lemma}
  Let $C_A$ be some $\gg$-valued multilinear form that transforms
  according to 
\[
  C_{u^*A}=u^{-1}C_Au.
\]
Then we also have
\[
  D_{u^*A}\,C_{u^*A}=u^{-1}D_AC_Au.
\]
\end{lemma}

We omit the proof because it is almost literally the same as the
preceding one if one replaces $d$ by $D$.\qed

Since the curvature of a connection $A$ transforms like
$F_{u^*A}=u^{-1}F_A u$, the three preceding lemmas yield the 

\begin{corollary}
  For the curvature $F_A$ of a local connection $A$ of class
  $W^{n-1,2}$ and a $W^{2,n}$-gauge transformation $u$ we have
  \begin{equation*}
    d^{*\wedge(n-2)}_{u^*A}F_{u^*A}=u^{-1}(d_A^{*\wedge(n-2)}F_A) u
  \end{equation*}
  and
  \begin{equation*}
    D^{n-2}_{u^*A}F_{u^*A}=u^{-1}(D^{n-2}_AF_A) u.
  \end{equation*}
  In particular, this implies the gauge invariance of the functionals
  $Y_n$ and $Z_n$ in the form $Y_n(u^*A)=Y_n(A)$, respectively
  $Z_n(u^*A)=Z_n(A)$. 
\end{corollary}

\begin{remark} \rm The reader may have expected $W^{n,2}$ gauge transformations 
instead of $W^{2,n}$ in the corollary. However, by the preceding
lemmas, we see that $d^{*\wedge(n-2)}_{u^*A}F_{u^*A}$ and $D^{n-2}_{u^*A}F_{u^*A}$
are defined even if $u$ is only in $W^{2,n}$. This is seen by iterating
arguments like $D_{u^*A}F_{u^*A}=D_{u^*A}(u^{-1}F_Au)$, where only one derivative
of $u$ is needed on the right-hand side.
\end{remark}

We have remarked in the introduction that $d_Ad_A$ and $d_A^*d_A^*$, even
when applied to forms like $A$ or $F_A$, are differential operators
of order $0$. More precisely, we have

\begin{lemma}\la{dd}
For all $\gg$-valued $2$-forms $B$ and $\gg$-valued $1$-forms $C$, we have 
the identities
\begin{align*}
  d_A^*d_A^*B&=-\st [F_A,\st B],\\
  d_Ad_AC&=[F_A,C].
\end{align*}
\end{lemma}

{\bf Proof.} The second assertion is more or less the definition of $F_A$,
and moreover a simpler variant of the proof of the first assertion, which
we now give.

We use $\st \st =(-1)^{k(m+1)}$ and $d^*=(-1)^{(k+1)m+1}
\st d\st $ when operating on $k$-forms, and $[X,Y]=(-1)^{k\ell+1}[Y,X]$ when $X$ is
a $k$-form and $Y$ is an $\ell$-form. Therefore
\eas
  d_A^*d_A^*B
  &=&-(-1)^md_A^*(\st d\st B+\st [A,\st B])\\
  &=&(-1)^m\Big(\st d\st \st d\st B+\st d\st \st [A,\st B]\\
  &&\qquad\quad{}+\st [A,\st \st d\st B]
    +\st [A,\st \st [A,\st B]]\Big)\\
  &=&-\st \Big(dd\st B+d[A,\st B]+[A,d\st B]+[A,[A,\st B]]\Big)\\
  &=&-\st \Big([dA,\st B]-[A,d\st B]+[A,d\st B]
    +[A,[A,\st B]]\Big)\\
  &=&-\st \Big([dA,\st B]+[A,[A,\st B]]\Big).
\eeas
The Jacobi identity (with correct signs) yields
\[
  [A,[A,\st B]]+(-1)^{m-1}[A,[\st B,A]]+[\st B,[A,A]]=0,
\]
from which we infer
\[
  2[A,[A,\st B]]=-[\st B,[A,A]]=[[A,A],\st B].
\]
We insert this in our previous calculation to find
\[
  d_A^*d_A^*B
  =-*\Big([dA,\st B]+\frac12\,[[A,A],\st B]\Big)
  =-\st [F_A,\st B]
\]
as asserted.\qed

\subsection{Gagliardo-Nirenberg interpolation}

In order to deal with lower order derivatives, we rely on the
Gagliardo-Nirenberg interpolation inequality in the following form
(see \cite[Thm. 1]{Ni}). 

\begin{theorem}\label{thm:interpol}
  Let $\Omega\subset\R^m$ be a bounded domain with the cone property and $u\in
  W^{k.r}(\Omega)\cap L^q(\Omega)$, where $1\le p,q\le\infty$ and $k\in\N$. Then we have the inequality
  \begin{equation*}
    \|D^ju\|_{L^p(\Omega)}\le C\|D^ku\|_{L^r(\Omega)}^{j/k}\|u\|_{L^q(\Omega)}^{1-j/k}+C\|u\|_{L^q(\Omega)}
  \end{equation*}
  provided $0\le j<k$ and $\frac1p=\frac jk\frac1r+(1-\frac
  jk)\frac1q$. Here, the constant $C$ depends only on $\Omega,k,j,q,$
  and $r$. 
\end{theorem}

\section{Gauge invariant interpolation and Sobolev inequalities} 
\label{sec:WA}

We define the Sobolev spaces $W^{k,p}_A(M,\wedge^2(T^*M)\otimes\gg_P)$
containing all sections $B$ of the bundle for which
\[
  \|B\|_{W^{k,p}_A(M)}:=\sum_{j=0}^k\|D_A^jB\|_{L^p(M)}
\]
is finite. The main purpose of this section is to verify that both
$Y_n(A)$ and $Z_n(A)$ can control $\|F_A\|_{W^{n-2,2}_A(M)}$.

\begin{theorem}\la{WA-estimates}
(i) For any smooth section $A$ of $T^*M\otimes\gg_p$ for which $Y_n(A)$
is finite, 
we have $F_A\in W^{n-2,2}_A(M,\wedge^2(T^*M)\otimes\gg_P)$, and even
\[
  \sum_{k=0}^{n-2}\|D_A^kF_A\|_{L^{2n/(k+2)}}^{2n/(k+2)}\le C(Y_n(A)+Y_n(A)^{2/n}).
\]
(ii) Assume that $m=\dim M\le 2n$.
For any smooth section $A$ of $T^*M\otimes\gg_p$ for which $Z_n(A)$
is finite, 
we have $F_A\in W^{n-2,2}_A(M,\wedge^2(T^*M)\otimes\gg_P)$, and even
\[
  \sum_{k=0}^{n-2}\|D_A^kF_A\|_{L^{2n/(k+2)}}^{2n/(k+2)}\le C(Z_n(A)+Z_n(A)^{n/2}).
\]
In both cases, the constant $C$ depends on $M$ and $n$ only.
\end{theorem}

{\bf Proof of (i).}
We remind the reader
of the Weitzenb\"ock formula, which somewhat symbolically reads
\[
  (d_A^*d_A+d_Ad_A^*-D_A^*D_A)B=F_A\#B+R_M\#B,
\]
where $F_A\#B$ and $R_M\#B$ mean universal bilinear forms applied to
$F_A$ or $R_M$, and $B$. Here $R_M$ is the Riemannian curvature of $M$.
Applying this to $B=F_A$ and multiplying with $F_A$, we find, using also
Bianchi's identity $d_AF_A=0$,
\eas
  \|D_AF_A\|_{L^2}^2&=&\int_M\<F_A,D_A^*D_AF_A\>\,dx\\
  &=&\int_M\<F_A,d_Ad_A^*F_A+F_A\#F_A+R_M\#F_A\>\,dx\\
  &\le&\|d_A^*F_A\|_{L^2}^2+\|F_A\#F_A\#F_A+R_M\#F_A\#F_A\|_{L^1},
\eeas
which we write as
\[
  \|D_AF_A\|_{L^2}^2\le\|d_A^*F_A\|_{L^2}^2
    +C\|F_A\#F_A\#F_A+F_A\#F_A\|_{L^1},
\]
since $R_M$ is a given smooth form. Now we proceed doing a similar calculation
for $D_A^2F_A$. We use, in that order (a) Weitzenb\"ock with $B=D_AF_A$, 
(b) $d_AD_AB-D_Ad_AB=F_A\#B$, $d_A^*D_AB-D_Ad_A^*B=F_A\#B$ and $d_AF_A=0$ by 
Bianchi, (c) Weitzenb\"ock with $B=d_A^*F_A$, (d) $d_A^*d_A^*B=F_A\#B$ using
Lemma \ref{dd} and the estimate \re{int0}. Those give
\eas
  \lll{\|D_A^2F_A\|_{L^2}^2}\\
  &\le&\|d_A^*D_AF_A\|_{L^2}^2+\|d_AD_AF_A\|_{L^2}^2
    +C\|D_AF_A\#D_AF_A\#F_A+D_AF_A\#D_AF_A\|_{L^1}\\
  &\le&C\|D_Ad_A^*F_A\|_{L^2}^2+C\|F_A\#F_A\|_{L^2}^2
    +C\|D_AF_A\#D_AF_A\#F_A+D_AF_A\#D_AF_A\|_{L^1}\\
  &\le&C\|d_Ad_A^*F_A\|_{L_2}^2+C\|d_A^*d_A^*F_A\|_{L^2}^2
    +C\|D_AF_A\#D_AF_A\#F_A+D_AF_A\#D_AF_A\|_{L^1}\\
  &&{}+C\|F_A\#F_A\|_{L^2}^2+C\|d_A^*F_A\#d_A^*F_A\#F_A+d_A^*F_A\#d_A^*F_A\|_{L^1}\\
  &\le&C\|d_Ad_A^*F_A\|_{L^2}^2+C\|LOT\|_{L^1}
\eeas
where the lower order terms $LOT$ are of the types
\eas
  &&F_A\#F_A\#F_A\#F_A,\quad F_A\#F_A\#F_A,\quad F_A\#F_A,\\
  &&D_AF_A\#D_AF_A\#F_A,\quad D_AF_A\#D_AF_A,
    \quad D_AF_A\#F_A\#F_A,\quad D_AF_A\#F_A.
\eeas 
Iterating this by a long but straightforward induction argument, we
find
\[
  \|D_A^{k-2}F_A\|_{L^2}^2\le C\|\das{k-2}F_A\|_{L^2}^2+C\|LOT\|_{L^1}
\]
for every $k\in\{2,\ldots,n\}$. Here the lower order terms each are of the 
form
\[
  \Kreuz_{j=0}^{k-3}(D^j_AF_A)^{\#h_j}=:TLOT
\]
(meaning ``this lower order term'') with
\[
  4\le s:=\sum_{j=0}^{k-3}(j+2)h_j\le2k,\qquad\sum_{j=0}^{k-3}h_j\ge2.
\]
Now we distinguish three cases.

{\bf Case 1.} {\em Only one of the $h_j$ is $\ne0$. }Then 
$2\le h_j\le\frac{2k}{j+2}$ and it is elementary to estimate
\[
  \|TLOT\|_{L^1}\le C\|D_A^jF_A\|_{L^{h_j}}^{h_j}\le C\Big(\|D_A^jF_A\|_{L^2}^2+
    \|D_A^jF_A\|_{L^{2k/(j+2)}}^{2k/(j+2)}\Big).
\]

{\bf Case 2.} {\em For all $j$ with $h_j\ne0$, we have $j+2\le\frac{s}2$.
}Then we use Young's inequality with exponents $\frac{s}{(j+2)h_j}$ whenever
$h_j\ne0$, and find, using $s\le2k$ and $\frac{s}{j+2}\ge2$,
\[
  \|TLOT\|_{L^1}\le C\sum_{j=0}^{k-3}\|D_A^jF_A\|_{L^{s/(j+2)}}^{s/(j+2)}
    \le C\sum_{j=0}^{k-3}\Big(\|D_A^jF_A\|_{L^2}^2+\|D_A^jF_A\|_{L^{2k/(j+2)}}^{2k/(j+2)}
    \Big).
\]

{\bf Case 3.} {\em The largest $j$ with $h_j\ne0$, let us call it $J$,
satisfies $J+2>\frac{s}2$. }Then $j+2\le\frac{s}2$ for all other
terms with $h_j\ne0$, and we must have $h_J=1$. We use Young's
inequality with exponents $2$ for the $J$-term, and $\frac{2(s-J-2)}{(j+2)h_j}$
for the others. The exponents $\frac{2(s-J-2)}{j+2}$ that occur in
the following calculation are $\ge 2$ because $j+2\le s-J-2$, and they are
$\le\frac{2k}{j+2}$ because $s-J-2<s-\frac{s}2\le k$. This justifies the
estimate
\eas
  \|TLOT\|_{L^1}&\le&C\|D_A^JF_A\|_{L^2}^2
    +C\sum_{j<J,\,h_j\ne0}\|D_A^jF_A\|_{L^{2(s-J-2)/(j+2)}}^{2(s-J-2)/(j+2)}\\
  &\le&C\sum_{j=0}^{k-3}\Big(\|D_A^jF_A\|_{L^2}^2+\|D_A^jF_A\|_{L^{2k/(j+2)}}^{2k/(j+2)}
    \Big).
\eeas

Summing up over all lower order terms therefore gives
\eq{int1}
  \|D_A^{k-2}F_A\|_{L^2}^2\le\|\das{k-2}F_A\|_{L^2}^2
    +C\sum_{j=0}^{k-3}\Big(\|D_A^jF_A\|_{L^2}^2+\|D_A^jF_A\|_{L^{2k/(j+2)}}^{2k/(j+2)}
    \Big).
\eeq

Now we start our interpolation considerations by remarking that
\eq{int2}
  \|\das{k}F_A\|_{L^2}\le\eps\|\das{n-2}F_A\|_{L^2}+C(\eps)\|F_A\|_{L^2}
\eeq
for $k\in\{0,\ldots,n-3\}$ is straightforward. Similarly,
\eq{int2a}
  \|D_A^kF_A\|_{L^2}\le\eps\|D_A^{n-2}F_A\|_{L^2}+C(\eps)\|F_A\|_{L^2},
\eeq
which is still easy to prove, but slightly more subtle.
This is in fact an easier variant of the following one. 
It reads
\eas
  \lll{\|D_A^kF_A\|_{L^{2n/(k+2)}}^{2n/(k+2)}
  =\int_M\<D_A^kF_A,|D_A^kF_A|^{\frac{2n}{k+2}-2}D_A^kF_A\>\,dx}\\
  &=&\int_M\<D_A^{k-1}F_A,|D_A^kF_A|^{\frac{2n}{k+2}-2}D_A^*D_A^kF_A\>\,dx\\
  &&{}+C\int_M\Big\<D_A^{k-1}F_A,|D_A^kF_A|^{\frac{2n}{k+2}-4}
    \<\<D_A^{k+1}F_A,D_A^kF_A\>,D_A^kF_A\>\Big\>\,dx\\
  &\le&C(\eps)\|D_A^{k-1}F_A\|_{L^{2n/(k+1)}}^{2n/(k+1)}
    +\eps\|D_A^{k+1}F_A\|_{L^{2n/(k+3)}}^{2n/(k+3)}\\
  &&{}+\eps\Big(\|D_A^kF_A\|_{L^{2n/(k+2)}}^{2n/(k+2)}
    +\|D_A^kF_A\|_{L^{2n/(k+3)}}^{2n/(k+3)}\Big)\\
  &\le&C(\eps)\|D_A^{k-1}F_A\|_{L^{2n/(k+1)}}^{2n/(k+1)}
    +\eps\|D_A^{k+1}F_A\|_{L^{2n/(k+3)}}^{2n/(k+3)}\\
  &&{}+\eps\Big(\|D_A^kF_A\|_{L^{2n/(k+2)}}^{2n/(k+2)}
    +\|D_A^kF_A\|_{L^2}^2\Big)
\eeas
At the second ``$=$'', 
we have used that $D_A$ is a metric connection. At exactly that point,
the version of it using $\das{k}$ instead of $D_A^k$ would fail, since we
would get one $D_A$, anyway. In some sense, this is why the proof of (i)
is not straightforward. For the first ``$\le$'', we have used Young's
inequality and \re{int0}.

Absorbing the second-last term into the left-hand side and estimating the
last one with \re{int2a}, we iterate and find
\eq{int3}
  \|D_A^kF_A\|_{L^{2n/(k+2)}}^{2n/(k+2)}\le\eps\|D_A^{n-2}F_A\|_{L^2}^2
    +C(\eps)(\|F_A\|_{L^2}^2+\|F_A\|_{L^n}^n)
\eeq
for all $k\in\{0,\ldots,n-3\}$.
Using \re{int3}, \re{int1}, \re{int3} again, and then 
\re{int2a}, 
we have
\eas
  \lll{\sum_{k=0}^{n-2}\|D_A^kF_A\|_{L^{2n/(k+2)}}^{2n/(k+2)}
  \le C(\|D_A^{n-2}F_A\|_{L^2}^2+\|F_A\|_{L^2}^2)+C\|F_A\|_{L^n}^n}\\
  &\le&C\Big(\|\das{n-2}F_A\|_{L^2}^2
    +\sum_{k=0}^{n-3}(\|D_A^kF_A\|_{L^2}^2
    +\|D_A^kF_A\|_{L^{2n/(k+2)}}^{2n/(k+2)})\Big)\\
  &\le&C(\eps)\Big(\|\das{n-2}F_A\|_{L^2}^2
    +\sum_{k=0}^{n-3}\|D_A^kF_A\|_{L^2}^2+\|F_A\|_{L^n}^n\Big)
    +\eps\|D_A^{n-2}F_A\|_{L_2}^2\\
  &\le&C(\eps)(\|\das{n-2}F_A\|_{L^2}^2+\|F_A\|_{L^2}^2
    +\|F_A\|_{L^n}^n)+C\eps\|D_A^{n-2}F_A\|_{L_2}^2.
\eeas
Absorbing the last term, we have proven assertion (i).\\[-2mm]\strut

{\bf Proof of (ii).} This is easier, and interesting in its own right. 
It is well-known that Sobolev inequalities for $D_A$ hold with constants 
not depending on $A$. This is because
\[
  2\,|B|\,|D|B||=|D|B|^2|=|D\<B,B\>|=2\,|\<B,D_AB\>|\le2\,|B|\,|D_AB|    
\]
implies $|D|B||\le|D_AB|$ for all $A$, and hence
\[
  \|B\|_{L^{p^*}}=\|\,|B|\,\|_{L^{p^*}}\le C_S\|\,D|B|\,\|_{W^{1,p}}
  \le C_S\|D_AB\|_{W^{1,p}}.
\]
Iterating, we also have higher order Sobolev inequalities with constants
independent of $A$, and in particular, using $\dim M\le 2n$,
\eas
  \sum_{k=0}^{n-2}\|D_A^kF_A\|_{L^{2n/(k+2)}}&\le& C\|F_A\|_{W^{n-2,2}_A}
  \;\le\; C(\|D^{n-2}_AF_A\|_{L^2}+\|F_A\|_{L^2})\\
  &\le& C Z_n(A)^{1/2},
\eeas
where the second ``$\le$'' comes from \re{int2a}.
This completes the proof of Theorem \ref{WA-estimates}.\qed

\begin{remark}\la{WA-local}
On manifolds with boundary, in particular on balls in $M$, we get similar
assertions. We have to use cutoff functions near the boundary, however, and
therefore get only $F_A\in W^{n-2,2}_A$ locally, with estimates on every
compact subset away from the boundary.
\end{remark}

\section{Uhlenbeck type estimates}\label{sec:uhlenbeck}




In dimension $m=2n$, Uhlenbeck \cite{Uh2}
showed that smallness of $\|F_A\|_{L^n(B_r)}$
ensures that by a suitable gauge change,
we can achieve $d^*A=0$,  $*A|_{\partial B_r}=0$ and 
\[
  \|A\|_{L^{2n}(B_r)}+\|DA\|_{L^n(B_r)}\le c\|F_A\|_{L^n(B_r)}.
\]

We will need something similar for higher order.

\begin{theorem}[higher order Uhlenbeck estimates]\label{uhl_n}
Assume $n\ge2$, $m\le 2n$, and let $B_r=B_r(a)$ be any ball of radius 
$r\in(0,1)$
in $M$. There is a constant $\kappa=\kappa(M)\in(0,1)$ such that the
following holds. Assume that $A\in W^{1,n}(B_r)$ satisfies
$A\in W^{n-1,2}_\loc(B_r\setminus\{0\})$ and
\begin{equation}
  \label{quasi-Sobolev}
  D^\ell_AF_A\in L^{2n/(\ell+2)}(B_r)
  \quad\mbox{for }\ell=1,\ldots,n-2
\end{equation}
(which is in particular satisfied in the case $A\in W^{n-1,2}(B_r)$).
Moreover, we assume that the curvature is small in the sense
\begin{equation}\label{smallness}
  r^{\frac{2n-m}{n}}\|F_A\|_{L^{n}(B_r)}< \kappa.
\end{equation}
Then the Uhlenbeck gauged version
$\Om$ of $A$ obeys the estimate
\begin{align}\label{Uhlenbeck-Higher-Est}
  &\sum_{\ell=0}^{n-1}r^{\frac{2n-m}{2n}(\ell+1)}\|D^\ell\Om\|_{L^{2n/(\ell+1)}(B_{r/2})}\\\nn
  &\qquad\le C r^{\frac{2n-m}{2}}\|D_\Omega^{n-2} F_\Om\|_{L^2(B_r)}
   +Cr^{\frac{2n-m}{n}}\|F_\Om\|_{L^n(B_r)}
\end{align}
with a constant $C$ depending on $M$ only. (Note that the powers
of $r$ all disappear in the critical dimension $m=2n$.)
\end{theorem}

\begin{remark}\rm
The assumption $A\in W^{n-1,2}_\loc(B_r\setminus\{0\})$ is technical. It could be
replaced by any weaker assumption that ensures all terms in \re{quasi-Sobolev} 
to be well-defined a.e.
\end{remark}

{\bf Proof of the theorem.} 
By scaling invariance, we may restrict ourselves to the
case $r=1$. We trivialize the bundle over $B_1$, and we consider
$B_1$ to be the Euclidean ball $B_1=B_1(0)\subset\R^m$, equipped with
some Riemannian metric $\gamma$. All constants involving $\gamma$
can be chosen independent of the choice of the ball, because
$M$ is compact. 

We will later choose $\kappa>0$ not larger than in Uhlenbeck's gauge Theorem~\cite{Uh2}. 
Then our assumption $\|F_A\|_{L^n}<\kappa$ allows us to find a gauge
$u\in W^{2,n}(B_1,G)$ such that $\Omega:=u^*A\in W^{1,n}(B_1)$ satisfies
$d^*\Omega=0$ and moreover,
\begin{equation}
  \label{Uhlenbeck-bound}
  \|\Omega\|_{L^{2n}(B_1)}+\|D\Omega\|_{L^n(B_1)}\le C\|F_A\|_{L^n(B_1)}<C\kappa.
\end{equation}
From gauge invariance of the total derivatives $D^\ell_AF_A$ and
\eqref{quasi-Sobolev} for $r=1$, we know
\begin{equation*}
  \|D^\ell_\Omega F_\Omega\|_{L^{2n/(\ell+2)}(B_1)}
  =
  \|D^\ell_A F_A\|_{L^{2n/(\ell+2)}(B_1)}
  <\infty
\end{equation*}
for $\ell=0,\ldots,n-2$. 
The remainder of the proof is divided into three steps.\\[-1ex]

\textbf{Step 1: Controlling $D^kd\Omega$ by $D_\Omega^kF_\Omega$} for
$k=1,\ldots,n-2$.

  Writing $D$ for the total derivative, applied separately to the
  coefficients of $F_\Omega$, we have a relation of the form
  \begin{align*}
    DF_\Omega
    &= D_\Omega F_\Omega +F_\Omega\#\Omega+F_\Omega\#\Gamma\\
    &= D_\Omega F_\Omega +D\Omega\#\Omega+\Omega\#\Omega\#\Omega+D\Omega\#\Gamma+\Omega\#\Omega\#\Gamma,
  \end{align*}
  where $\Gamma$ represents the Christoffel symbols of the manifold $M$. 
  Keeping in mind that $\Gamma$ and all its derivatives are bounded by
  constants depending only on $M$, 
  we can generalize the last formula
  inductively to higher order, with the result
\begin{equation}
  \label{DF-Comparison}
  D^kF_\Omega 
  = 
  D_\Omega^kF_\Omega
  +
  \sum_{j\in J_{k}}\Kreuz_{i=1}^\ell D^{j_i-1}\Omega,
\end{equation}
for $k=1,\ldots,n-2$, where we abbreviated
\[
  J_k:=\{j=(j_1,\ldots,j_\ell):\,\ell\ge1,\,1\le j_i\le k+1,\,
    2\le j_1+\ldots+j_\ell\le  k+2\}.
\]
Applying first H\"older's inequality with exponents $(k+2)/j_i$ and
then Young's inequality with exponents $(j_1+\ldots+j_\ell)/j_i$, we infer
\begin{align*}
  &\|D^k F_\Omega\|_{L^{2n/(k+2)}(B_{\rho})}\\
  &\quad\le 
  C\|D_\Omega^k F_\Omega\|_{L^{2n/(k+2)}(B_{\rho})}
  +
  C\sum_{j\in
    J_k}\prod_{i=1}^\ell\|D^{j_i-1}\Omega\|_{L^{2n/j_i}(B_{\rho})}\\
  &\quad\le C\|D_\Omega^k F_\Omega\|_{L^{2n/(k+2)}(B_{\rho})}
  +
  C\sum_{\ell=0}^k\big(\|D^\ell\Omega\|_{L^{2n/(\ell+1)}(B_{\rho})}
    +\|D^\ell\Omega\|_{L^{2n/(\ell+1)}(B_{\rho})}^{(k+2)/(\ell+1)}\big)\,,
\end{align*}
for every $\rho\in[\frac12,1]$. 
From Leibnitz' rule and Young's inequality, we deduce 
\begin{equation*}
  \|D^k (\Omega\wedge\Omega)\|_{L^{2n/(k+2)}(B_{\rho})}
  \le
  C\sum_{\ell=0}^k \|D^\ell\Omega\|_{L^{2n/(\ell+1)}(B_{\rho})}^{(k+2)/(\ell+1)}\,.
\end{equation*}
Joining the two preceding estimates and keeping in mind
$F_\Omega=d\Omega+\Omega\wedge\Omega$, we arrive at
\begin{align}\label{dOmega-bound}
  &\|D^k d\Omega\|_{L^{2n/(k+2)}(B_\rho)}\\\nn
  &\quad\le 
  C\|D_\Omega^k F_\Omega\|_{L^{2n/(k+2)}(B_\rho)}
  +
  C\sum_{\ell=0}^k\big(\|D^\ell\Omega\|_{L^{2n/(\ell+1)}(B_{\rho})}
    +\|D^\ell\Omega\|_{L^{2n/(\ell+1)}(B_{\rho})}^{(k+2)/(\ell+1)}\big)
\end{align}
for $k=1,\ldots,n-2$.\\

\textbf{Step 2: Proof of $\Omega\in W^{n-1,2}_\loc(B_1)$.}

We will prove by induction over $k=1,\ldots,n-1$ that 
\begin{equation}
D^\ell\Omega\in L^{2n/(\ell+1)}_\loc(B_1)
\qquad\mbox{for }\ell=0,\ldots,k.\label{Sobolev-claim}
\end{equation}
For $k=1$, this is a consequence of
Uhlenbeck's result~(\ref{Uhlenbeck-bound}). Next, we assume
that~(\ref{Sobolev-claim}) is valid for $k\in\{1,\ldots,n-2\}$ and
wish to prove it for $k+1$. To this end, we calculate
\begin{equation*}
  -\Delta\Omega=d^*d\Omega+dd^*\Omega=d^*d\Omega,
\end{equation*}
which implies
\begin{align*}
  \|\Delta D^{k-1}\Omega\|_{L^{2n/(k+2)}(B_\rho)}
  &\le
  C\|d\Omega\|_{W^{k,2n/(k+2)}(B_\rho)}\\
  &\le
  C\|D^kd\Omega\|_{L^{2n/(k+2)}(B_\rho)}
  +
  C\|\Omega\|_{W^{k,2n/(k+2)}(B_\rho)}
  <\infty
\end{align*}
for every $\rho\in[\frac12,1]$. The finiteness of the right-hand side is a
consequence of~(\ref{dOmega-bound}) and the induction
assumption~(\ref{Sobolev-claim}). Now classical Calder\'on-Zygmund
theory implies $D^{k+1}\Omega\in L^{2n/(k+2)}_\loc(B_1)$. Proceeding in
this manner up to the order $k=n-1$, we arrive at $D^{n-1}\Omega\in
L^2_\loc(B_1)$, which establishes the claim $\Omega\in W^{n-1,2}_\loc(B_1)$.\\

\textbf{Step 3: Proof of estimate~(\ref{Uhlenbeck-Higher-Est}).}    

The Sobolev regularity established in the preceding step now justifies
the following calculations that will lead to the desired estimate. 
For given radii $R,S$ with $\frac12\le R<S\le \frac34$, we choose a cut-off function $\p\in
C^\infty_{\rm cpt}(B_S,[0,1])$ with $\p\equiv1$ on $B_R$ and
$\|D^j\p\|_{L^\infty}\le C/(S-R)^j$ for all
$j\in\{1,\ldots,n-1\}$. Using once more $d^*\Omega=0$, we obtain
\begin{equation*}
  \Delta(\p \Omega)=d^*\big(d\p\wedge\Omega +\p d\Omega \big)
  +d(d\p\cdot \Omega).
\end{equation*}
Differentiating this identity $(n-3)$ times, using the properties
of $\p$ and $S-R<1$, we deduce
\begin{equation*}
  \|\Delta D^{n-3}(\p \Omega)\|_{L^2}
  \le
  C\|D^{n-2}d\Omega\|_{L^2(B_S)}
  +
  C\sum_{\ell=0}^{n-2}\tfrac 1{(S-R)^{n-\ell-1}}\,\|D^\ell\Omega\|_{L^2(B_S)}\,.
\end{equation*}
Next, we apply (\ref{dOmega-bound}) for $k=n-2$ and $\rho=S$, with the result
\begin{align*}
  &\|\Delta D^{n-3}(\p \Omega)\|_{L^2}
  \le
  C\|D_\Omega^{n-2} F_\Omega\|_{L^2(B_S)}\\
  &\qquad\qquad\quad+
  C\sum_{\ell=0}^{n-2}\Big(
     \tfrac 1{(S-R)^{n-\ell-1}}\,\|D^\ell\Omega\|_{L^{2n/(\ell+1)}(B_S)}
    +\|D^\ell\Omega\|_{L^{2n/(\ell+1)}(B_S)}^{n/(\ell+1)}\Big). 
\end{align*}
At this stage, we once more apply classical Calder\'on-Zygmund theory,
which yields the bound
\begin{align}\label{Higher-Order-Bound}
  &\|D^{n-1}\Omega\|_{L^2(B_R)}
  \le
  C\|D_\Omega^{n-2} F_\Omega\|_{L^2(B_S)}\\\nn
  &\qquad\qquad\quad+
  C\sum_{\ell=0}^{n-2}\Big(
     \tfrac 1{(S-R)^{n-\ell-1}}\,\|D^\ell\Omega\|_{L^{2n/(\ell+1)}(B_S)}
    +\|D^\ell\Omega\|_{L^{2n/(\ell+1)}(B_S)}^{n/(\ell+1)}\Big).
\end{align}
In order to bound the terms of the last sum, 
we apply Gagliardo-Nirenberg interpolation in the
form stated in Theorem~\ref{thm:interpol}, which gives
\begin{align*}
  \|D^\ell\Omega\|_{L^{2n/(\ell+1)}(B_S)}
  &\le
  C\Big(\|D^{n-1}\Omega\|_{L^2(B_S)}^{\ell/(n-1)}\|\Omega\|_{L^{2n}(B_S)}^{(n-\ell-1)/(n-1)}
 +\|\Omega\|_{L^{2n}(B_S)}\Big)
\end{align*}
for every $\ell\in\{0,\ldots,n-2\}$.
For a parameter $\eps>0$ to be chosen later, we twice apply
Young's inequality to the right-hand side, which yields 
\begin{equation*}
  \tfrac 1{(S-R)^{n-\ell-1}}\|D^\ell\Omega\|_{L^{2n/(\ell+1)}(B_S)}
  \le
  \eps \|D^{n-1}\Omega\|_{L^2(B_S)}
  +
  \frac{C(\eps)}{(S-R)^{n-1}}\,\|\Omega\|_{L^{2n}(B_S)},
\end{equation*}
as well as 
\begin{equation*}
  \|D^\ell\Omega\|_{L^{2n/(\ell+1)}(B_S)}^{n/(\ell+1)}
  \le
  \eps \|D^{n-1}\Omega\|_{L^2(B_S)}
  +C(\eps)\|\Omega\|_{L^{2n}(B_S)}^n.
\end{equation*}
We recall that $\|\Omega\|_{L^{2n}(B_S)}\le C\kappa$ with
$\kappa<1$, so that we also may drop the exponent $n$ in the last
term. Plugging the two preceding estimates
into~(\ref{Higher-Order-Bound}), we obtain the bound  
\begin{align*}
  &\|D^{n-1}\Omega\|_{L^2(B_R)}\\
  &\qquad\le
  C\eps \|D^{n-1}\Omega\|_{L^2(B_S)}
  +
   C\|D_\Omega^{n-2} F_\Omega\|_{L^2(B_1)}
   +
   \frac{C(\eps)}{(S-R)^{n-1}}\,\|\Omega\|_{L^{2n}(B_1)}\,.
\end{align*}
Now we choose $\eps>0$ so small that $C\eps\le\frac12$.
Recalling that the above estimate holds
for any $\frac12\le R<S\le \frac34$, we iterate it in a standard way
(cf. \cite[Lemma V.3.1]{G}) to get
\begin{equation*}
  \|D^{n-1}\Omega\|_{L^2(B_{1/2})}
  \le
   C\|D_\Omega^{n-2} F_\Omega\|_{L^2(B_1)}
   +
   C\|\Omega\|_{L^{2n}(B_1)}\,.
\end{equation*}
Combining this with~(\ref{Uhlenbeck-bound}), we arrive at
\begin{equation*}
  \|D^{n-1}\Omega\|_{L^2(B_{1/2})}
  +
  \|\Omega\|_{L^{2n}(B_1)}
  \le
   C\|D_\Omega^{n-2} F_\Omega\|_{L^2(B_1)}
   +
   C\|F_\Omega\|_{L^n(B_1)}.
\end{equation*}
This implies the claimed estimate~(\ref{Uhlenbeck-Higher-Est}) by
another application of the Gagliardo-Nirenberg interpolation estimate.
\qed

\section{Removability of point singularities of the underlying bundle}
\label{sec:rem}
For Yang-Mills connections, there exist removable singularity results
in $4$ dimensions \cite{Uh1} and in higher dimensions \cite{TT}. 
A related partial regularity result for
Yang-Mills connections in higher dimensions has been established in \cite{MR}. In our
case, we do not work with Yang-Mills connections. Therefore it is important
to note that Petrache and Rivi\`ere have removed the assumption of
having a Yang-Mills connection from the removability result.

In this section, we prove a higher order version of their
ground-breaking result \cite[Thm. 3.2]{PR} for higher dimensions. The crucial
observation that enables Petrache and Rivi\`ere to construct a local 
trivialization without using the Yang-Mills equation 
is the continuity of gauge transformations between Coulomb gauges.
This follows from estimates involving Lorentz spaces, a technique that
has been introduced to Yang-Mills theory by Rivi\`ere in \cite{Ri}.
We generalize the arguments from \cite[Thm. 3.2]{PR}, which are the $n=2$ 
case of

\begin{proposition}\la{EichZwCoul}
Let $n\ge2$, $U\subset\R^{2n}$ be a bounded domain.
Assume that $A$ and $B=u^{-1}Au+u^{-1}du$
are ${\mathfrak g}$-valued $1$-forms of class $W^{n-1,2}$ satisfying
\[
  d^*A=d^*B=0 \qquad\mbox{on $U$. }
\]

(i) Then the gauge change $u$ is in $W^{n,2}_\loc\cap C^0_\loc(U,G)$, 
and for any $V\subset\!\!\!\subset U$, there is some constant $\bar{u}\in G$
depending on $u$, such that
\begin{align*}
  &\|u-\bar{u}\|_{W^{n,2}\cap C^0(V)}\\[1ex]
  &\qquad\le
  C(\|A\|_{W^{n-1,2}(U)}+\|B\|_{W^{n-1,2}(U)}+\|A\|_{W^{n-1,2}(U)}^n+\|B\|_{W^{n-1,2}(U)}^n),
\end{align*}
where $C$ depends only on $n$, $U$, $V$, and $G$.

(ii) Moreover, there is a $\delta>0$ depending on $n$, $U$, $V$, and
$G$ such that in the case 
\begin{equation}
  \|A\|_{W^{n-1,2}(U)}+\|B\|_{W^{n-1,2}(U)}<\delta,\label{small-norm-assumption}
\end{equation}
the following holds.
For any $W\subset\!\!\!\subset V$, there is another gauge change
$\tilde{u}\in W^{n,2}\cap C^0(U,G)$ that coincides with $u$ on $W$ and with the
constant $\bar{u}$ on $U\setminus V$, with the estimate
\begin{align*}
  &\|\tilde{u}-\bar{u}\|_{W^{n,2}\cap C^0(U)}\\[1ex]
  &\qquad\le
  C(\|A\|_{W^{n-1,2}(U)}+\|B\|_{W^{n-1,2}(U)}+\|A\|_{W^{n-1,2}(U)}^n+\|B\|_{W^{n-1,2}(U)}^n).
\end{align*}
The constant $C$ here depends additionally on $W$.
\end{proposition}

{\bf Proof of (i), }modelled after \cite[Prop. 3.5]{PR}.

If no domain is indicated in integral norms, it is assumed to be $U$.
By the Lorentz space version of the Sobolev embedding, we have
$A,B\in L^{(2n,2)}$ and
\[
  \|A\|_{L^{(2n,2)}}\le C\|A\|_{W^{n-1,2}},\qquad \|B\|_{L^{(2n,2)}}\le C\|B\|_{W^{n-1,2}}.
\]
Since
\[
  du=uB-Au,
\]
and since $u\in L^\infty$ as it takes values in $G$, the previous estimates
imply $du\in L^{(2n,2)}$ and
\[
  \|du\|_{L^{(2n,2)}}\le C(\|A\|_{L^{(2n,2)}}+\|B\|_{L^{(2n,2)}})
  \le C(\|A\|_{W^{n-1,2}}+\|B\|_{W^{n-1,2}}).
\]
Leibnitz' rule for functions $h$ and $1$-forms $\omega$ is
$d^*(h\omega)=hd^*\omega-dh\cdot\omega$, and this carries over to the 
${\mathfrak g}$-valued case. Using $d^*A=d^*B=0$, we therefore find
\[
  -\Delta u=d^*(du)=d^*(uB-Au)=du\cdot A-B\cdot du.
\]
On the right-hand side, we have products of $L^{(2n,2)}$-functions,
which are in $L^{(n,1)}$. Hence
\begin{align}\label{rem0}\nn
  \|\Delta u\|_{L^{(n,1)}}
  &\le\|du\|_{L^{(2n,2)}}(\|A\|_{L^{(2n,2)}}+\|B\|_{L^{(2n,2)}})\\
  &\le C(\|A\|_{W^{n-1,2}}^2+\|B\|_{W^{n-1,2}}^2).
\end{align}
Let $\eta$ be a $V$-$U$-cutoff function. Then, writing $u_0$ for
the mean value of $u$ in some $\R^{k\times k}\supset G$,
\begin{align}
  \label{rem1}
  \|\Delta(\eta(u-u_0))\|_{L^{(n,1)}}
  &\le\|\Delta u\|_{L^{(n,1)}}+C(\|du\|_{L^{(n,1)}}+\|u-u_0\|_{L^{(n,1)}})\\
  &\le\|\Delta u\|_{L^{(n,1)}}+C\|du\|_{L^{(2n,2)}}\nn\\
  &\le C(\|A\|_{W^{n-1,2}}+\|B\|_{W^{n-1,2}}+\|A\|_{W^{n-1,2}}^2+\|B\|_{W^{n-1,2}}^2)\nn\\
  &=:C(...),\nn
\end{align}
and $\eta(u-u_0)\equiv0$ near $\partial U$. The standard elliptic 
estimate in Lorentz space gives
\[
  \|u-u_0\|_{W^{2,(n,1)}(V)}\le C\|\eta(u-u_0)\|_{W^{2,(n,1)}(U)}\\
  \le C(...).
\]
Using $W^{2,(n,1)}\hookrightarrow L^\infty$ (which would not hold for $W^{2,n}$),
we also have
\[
  \|u-u_0\|_{L^\infty(V)}\le C(...).
\]
In order to show that $u$ is continuous, we 
use the estimate~\eqref{rem0} with $U$ replaced by an arbitrary 
ball $B_{2\rho}(x_0)\subset U$  and $V$ replaced by
$B_\rho(x_0)$. Abbreviating $u_\rho:=\mint_{B_{2\rho}(x_0)}u$ and bounding the lower order terms on the right-hand side
by H\"older's inequality, we deduce 
\begin{align*}
  \|u-u_\rho\|_{L^\infty(B_{\rho}(x_0))}
  &\le
  C\sum_{k=0}^{n-1}\Big(\|D^kA\|_{L^\frac{2n}{k+1}(B_{2\rho}(x_0))}
  +\|D^kB\|_{L^\frac{2n}{k+1}(B_{2\rho}(x_0))}\\
  &\phantom{C\sum_{k=0}^{n-1}\Big(}
  +\|D^kA\|^2_{L^\frac{2n}{k+1}(B_{2\rho}(x_0))}
  +\|D^kB\|^2_{L^\frac{2n}{k+1}(B_{2\rho}(x_0))}\Big)\,.
\end{align*}
In this form, both sides of the inequality are scaling invariant, from
which we infer that the constant $C$ can be chosen independently from
$\rho>0$. Keeping in mind the Sobolev embedding
$W^{n-1,2}\hookrightarrow W^{k,2n/(k+1)}$ that holds for any
$k=0,\ldots,n-1$, we deduce
$\|u-u_\rho\|_{L^\infty(B_{\rho}(x_0))}\to 0$ as $\rho\searrow0$,  
from which we see that $u$ is really continuous. 

Observe that $u\in G$ everywhere, and hence
\[
  \dist(u_0,G)\le C(...).
\]
This means that there is also $\bar{u}\in G$ such that
\[
  \|u-\bar{u}\|_{L^\infty(V)}\le C(...).
\]
Now that $u-\bar{u}$ is estimated in $W^{2,n}\cap L^{\infty}$, all that is missing
are higher order estimates for $u$. Starting with $u\in W^{2,n}\cap L^{\infty}$,
we plug that into $du=uB-Au$, and iterate that with any better result we
achieve that way, consecutively finding (with $V$ as domain of integration)
\eas
  du&\in&(W^{2,n}\cap L^\infty)\cdot W^{2,\frac{2n}3}\hookrightarrow W^{2,\frac{2n}3},\\
  du&\in&(W^{3,\frac{2n}3}\cap L^\infty)\cdot W^{3,\frac{2n}4}\hookrightarrow W^{3,\frac{2n}4},\\
  &\vdots&\\
  du&\in&(W^{n-1,\frac{2n}{n-1}}\cap L^\infty)\cdot W^{n-1,2}\hookrightarrow W^{n-1,2}.
\eeas
The last one implies $u\in W^{n,2}(V)$, with the asserted estimates, in which 
the powers of the norms build up because of iterated multiplication.

{\bf Proof of (ii)}. Now we assume that (\ref{small-norm-assumption}) holds 
for some sufficiently small $\del>0$. If we choose $\delta>0$ small
enough in dependence on $n,U,V$, and $G$, the estimate from (i) shows that 
the image of $\bar{u}^{-1}u_{|V}$ is contained in some neighborhood
of $e\in G$ on which $\exp_e^{-1}$ is defined and well behaved. More
precisely, we can assume that the first $n$ derivatives of $\exp_e^{-1}$ and $\exp_e$ on
this neighborhood or its exponential image are bounded by a constant depending
on $G$ only. Writing $h(x):=\exp_e^{-1}(\bar{u}^{-1}u(x))$, and using a 
$W$-$V$-cutoff function $\eta:U\to[0,1]$, we let
\[
  \tilde{u}(x):=\bar{u}\exp_e(\eta(x)h(x)),
\]
and easily see that it has the asserted properties.\qed

{\bf Remark.} $A,B\in W^{1,2n}$ would not have been enough to infer 
$u\in C^0$, because $u\in W^{2,(n,1)}\hookrightarrow L^\infty$ is crucial, and we
would only get $W^{2,(n,n/2)}$ instead. Therefore, we do need $A,B\in W^{n-1,2}$.
For $n=2$ (which is the case treated in \cite{PR}), both conditions coincide.

\strut

Now we follow \cite{PR} in proving the following removable singularity
theorem. 

\begin{theorem}\la{RemSing} 
Let  $P$ be a principal bundle over $B^{2n}\setminus\{0\}$.
Assume we are given a connection $D_0+A$ on $P$ of class
$W^{n-1,2}_\loc(B^{2n}\setminus\{0\})$, for which
\[
  \sum_{\ell=0}^{n-2}\|D^\ell_AF_A\|_{L^{2n/(\ell+2)}}<\infty.
\]
Then there exists a gauge of class $W^{n,2}_\loc$ in which the bundle extends to a smooth bundle over
$B^{2n}$, and the connection extends to a connection in $W^{n-1,2}(B^{2n})$
in this gauge.
\end{theorem}

For what follows, we define four sequences of spherical shells,
\eas
  Q_k&:=&B_{2^{-8k-3}}\setminus B_{2^{-8k-14}}\,,\\
  R_k&:=&B_{2^{-8k-2}}\setminus B_{2^{-8k-15}}\,,\\
  S_k&:=&B_{2^{-8k-1}}\setminus B_{2^{-8k-16}}\,,\\
  T_k&:=&B_{2^{-8k}}\setminus B_{2^{-8k-17}}\,
\eeas
for all $k\in\N\cup\{0\}$. Note $Q_k\subset R_k\subset S_k\subset T_k$.

\begin{lemma}\la{RingLemma} 
There exists a constant $\del\in(0,1)$ such that
\[
  \|D^{n-2}_AF_A\|_{L^2(T_k)}+\|F_A\|_{L^n(T_k)}<\del
\]
implies that the bundle $E$ is trivial over $Q_k$, and that in a 
suitable gauge the connection is represented by a $W^{n-1,2}$-form $A_k$
with $d^*A_k=0$ and the estimate
\eq{rem3}
  \sum_{\ell=0}^{n-1}\|D^\ell A_k\|_{L^{\frac{2n}{\ell+1}}(Q_k)}
    \le C(\|D^{n-2}_AF_A\|_{L^2(T_k)}+\|F_A\|_{L^n(T_k)}).
\eeq
Here, the constants $\delta$ and $C$ depend only on $M$ and $G$,
and in particular not on $k$. 
\end{lemma}

{\bf Proof,} following \cite[Lemma 3.6]{PR}.

By scaling invariance, it is enough to prove the lemma for $k=0$. We cover
$S_0$ by two charts $S_+$ and $S_-$ both diffeomorphic to $B^{2n}$, e.g.\
\[
  S_+:=\{x\in S_0:x_{2n}>-2^{-18}\},\qquad S_-:=\{x\in S_0:x_{2n}<2^{-18}\}.
\] 
On both $S_+$ and $S_-$, we can apply the higher order Uhlenbeck
estimates from Theorem \ref{uhl_n}, which clearly hold also on domains 
diffeomorphic to a ball, with a constant additionally depending on the
diffeomorphism. Because of the scaling invariance of the assertion 
\eqref{rem3}, we have to choose only two diffeomorphisms for $S_+$ and $S_-$
before applying Theorem \ref{uhl_n}, which means that the additional 
constant will depend only on $n$.
Hence we can assume that the connection is
represented by $A_+$ on $S_+$ and by $A_-$ on $S_-$ such that
\begin{align*}
  &\sum_{j=0}^{n-1}\Big(\|D^jA_+\|_{L^{\frac{2n}{j+1}}(S_+)}+\|D^jA_-\|_{L^{\frac{2n}{j+1}}(S_-)}\Big)\\
  &\qquad\le C(\|D^{n-2}_AF_A\|_{L^2(T_0)}+\|F_A\|_{L^n(T_0)})
  \le C\del.
\end{align*}
Now we let $R_+:=R_0\cap S_+$, $R_-:=R_0\cap S_-$, and
\[
  R_{++}:=\{x\in R_0:x_{2n}>-2^{-19}\},\qquad
  R_{--}:=\{x\in R_0:x_{2n}<2^{-19}\}.
\]
Note that in particular $R_{++}\subset R_+$ and $R_{--}\subset R_-$.
Proposition \ref{EichZwCoul} gives us $A_+=u^{-1}du+u^{-1}A_-u$ for some gauge
transformation $u$ controlled in $W^{n,2}_\loc\cap C^0_\loc(S_+\cap S_-)$. Part
(ii) provides us with a modification $\tilde{u}\in W^{n,2}\cap C^0(R_{--},G)$
coinciding with $u$ on $R_{++}\cap R_{--}$ and with some constant
$\bar{u}\in G$
on $R_{--}\setminus R_+$, such that
\eq{uklein}
  \|\tilde{u}-\bar{u}\|_{W^{n,2}\cap C^0(R_{--})}
  \le 
  C(\|D^{n-2}_AF_A\|_{L^2(T_0)}+\|F_A\|_{L^n(T_0)})
  \le
  C\del.
\eeq
Now here is a representative of the connection on all of $R_0$:
\eq{anull}
  \widetilde{A}_0(x):=\left\{\begin{array}{ll}
    A_+(x)&\mbox{ if }x_{2n}\ge0,\\
    \tilde{u}^*A_-(x)&\mbox{ if }x_{2n}<0.
  \end{array}\right.
\eeq
It is important to say how to interpret that assertion. Originally, the 
connection described by $A_+$ and $A_-$ is given on a bundle over $R_0$
which is glued together along $R_0\cap\{x_{2n}=0\}$ with the transition
function $\tilde{u}=u$. By the estimate \re{uklein}, $u$ takes its
values in a small ball around $\bar{u}$ in $G$ and therefore represents
the trivial homotopy class in $[R_0\cap\{x_{2n}=0\},G]\cong\pi_{2n-2}(G)$.
Since $R_0$ retracts to $S^{2n-1}$, the $G$-principal fibre bundles over
$R_0$ are classified as those over $S^{2n-1}$, that is by the element of
$\pi_{2n-2}(G)$ that the transition map, when restricted to an equator,
represents. See \cite[Section~4.4]{Na} for details on the classification.
Now $u$ represents the trivial class, hence the original bundle must
have been trivial, and \re{anull} expresses the connection in a 
trivialization of the bundle, where the transition map is the identity. 
The bounds from \re{uklein} of $\tilde{u}$, together with the estimates for
$A_+$ and $A_-$, show 
\eq{aklein}
  \sum_{j=0}^{n-1}\|D^j\widetilde{A}_0\|_{L^{\frac{2n}{j+1}}(R_0)}
    \le 
    C(\|D^{n-2}_AF_A\|_{L^2(T_0)}+\|F_A\|_{L^n(T_0)})
    \le 
    C\del.
\eeq
That is almost all we would require from $A_0$, but $\widetilde{A}_0$
is not yet in Coulomb gauge, $d^*\widetilde{A}_0$ may be $\ne0$.
Trying to gauge the connection again, we cannot
apply Uhlenbeck's gauge theorem directly, since we do not know if it still
holds on domains which are not diffeomorphic to $B^{2n}$.
We therefore modify $\widetilde{A}_0$ to find a connection on a ball, namely
$\widehat{A}_0$ on $B_{1/4}$ given by
\[
  \widehat{A}_0:=\eta\widetilde{A}_0,
\]
where here $\eta$ is a radial cutoff function which is $\equiv1$ on $Q_0$
and $\equiv0$ on $B_{1/4}\setminus R_0$ (where $\widetilde{A}$ is undefined).
Then, by \re{aklein}, we have
\[
  \|F_{{\widehat A}_0}\|_{L^n(B_{1/4})}
    \le C(\|F_A\|_{L^n(R_0)}
     +\|\widetilde{A}_0\|_{L^{n}(R_0)}
     +\|\widetilde{A}_0\|_{L^{2n}(R_0)}^2)\le C\del,
\]
and if $\del$ has been chosen small enough, we may apply Uhlenbeck's gauge
theorem to find a gauge-transformed version $A_0$ of $\widehat{A}_0$ with
$d^*A_0=0$. And since $\widehat{A}_0=\widetilde{A}_0$ on $Q_0$, $A_0$ represents
our originally given connection on $Q_0$. The asserted estimate for $A_0$
follows again from Theorem \ref{uhl_n}.\qed

The {\bf proof of Theorem \ref{RemSing}} again follows the arguments outlined
in \cite{PR}.

By restricting to a very small ball $B_r$ centered at $0$, we may assume that
the integrals of $|F_A|^n$ and $|D_A^{n-2}F_A|^2$ are as small as we want.
Both of them are scaling invariant, hence we may rescale to $B^{2n}$ and
assume
\[
  \|D^{n-2}_AF_A\|_{L^2(B^{2n})}+\|F_A\|_{L^n(B^{2n})}<\del,
\]
for a small $\del\in(0,1)$ yet to be chosen. Letting
\[
  \del_k:=\|D^{n-2}_AF_A\|_{L^2(T_k)}+\|F_A\|_{L^n(T_k)}\,,
\]
we have 
\begin{align}\label{sum_delta}
  \sum_{k=0}^\infty\del_k^n
  &\le
  C\sum_{k=0}^\infty\|D^{n-2}_AF_A\|_{L^2(T_k)}^2 
  +
  C\sum_{k=0}^\infty\|F_A\|_{L^n(T_k)}^n\\
  &\le
  C(\|D^{n-2}_AF_A\|_{L^2(B^{2n})}^2 +\|F_A\|_{L^n(B^{2n})}^n)\nn\\
  &\le
  C\del^2.\nn
\end{align}
If $\del$ is small enough, we can
apply Lemma \ref{RingLemma} to find forms $A_k$ representing our connection
on $Q_k$ and satisfying $d^*A_k=0$ as well as the estimates
\eq{aksmall}
  \sum_{j=0}^{n-1}\|D^jA_k\|_{L^{\frac{2n}{j+1}}(Q_k)}\le C\del_k.
\eeq
Let $W_k:=B_{2^{-8k-12}}\setminus B_{2^{-8k-13}}$ be a shell strictly inside
$U_k:=Q_k\cap Q_{k+1}=B_{2^{-8k-11}}\setminus B_{2^{-8k-14}}$. Using Proposition
\ref{EichZwCoul}\,(ii) (which \re{aksmall} allows us), we find mappings 
$\tilde{u}_k\in W^{n,2}\cap C^0(U_k,G)$ and constants $\bar{u}_k\in G$ such that
\eas
  &&A_{k+1}=\tilde{u}_k^*A_k\qquad\mbox{ on }W_k,\\
  &&\tilde{u}_k\equiv\bar{u}_k\qquad\mbox{ near }\d U_k,\\
  &&\|\tilde{u}_k-\bar{u}_k\|_{W^{n,2}\cap C^0(U_k)}\le C\del_k.
\eeas
Defining
\[
  \bar{v}_k:=\bar{u}_{k-1}^{-1}\ldots \bar{u}_1^{-1}\bar{u}_0^{-1},\qquad \Omega_k:=\bar{v}_k^*A_k,
\]
we have
\eq{udef}
  \Omega_{k+1}=u_k^*\Omega_k\mbox{ on }W_k,
   \mbox{ where }u_k:=\bar{v}_k^{-1}\tilde{u}_k\bar{u}_k^{-1}\bar{v}_k.
\eeq
The modified gauge changes $u_k$ satisfy
\[
  u_k-e=\bar{v}_k^{-1}(\tilde{u}_k-\bar{u}_k)\bar{u}_k^{-1}\bar{v}_k,
\]
and hence
\eq{uksmall}
  \|u_k-e\|_{W^{n,2}\cap C^0(U_k)}=\|\tilde{u}_k-\bar{u}_k\|_{W^{n,2}\cap C^0(U_k)}
  \le C\del_k\to0.
\eeq
Together with \re{aksmall}, this implies
\eq{oksmall}
  \sum_{j=0}^{n-1}\|D^j\Omega_k\|_{L^{\frac{2n}{j+1}}(U_k)}\le C\del_k.
\eeq
Not only are the $u_k$ close to the identity of $G$, but even $u_k\equiv e$
near $\d U_k$. Together with \re{udef}, this means that 
\[
  \Omega:=\left\{\begin{array}{ll}
    \Omega_k&\mbox{ on every }B_{2^{-8k-6}}\setminus B_{2^{-8k-11}}\,,\\
    u_k^*\Omega_k&\mbox{ on every }B_{2^{-8k-11}}\setminus B_{2^{-8k-14}}
  \end{array}\right.
\]
defines a connection form $\Omega$ on $B_{1/32}\setminus\{0\}$ which is 
locally (that is, away from $0$) in $W^{n-1,2}$. This is because $u_k$ has been
constructed carefully in order to prevent $\Omega$ from having jumps across
spheres. And $\Omega$ represents the originally given connection in some 
gauge, since it is obtained from the $A_k$ by a sequence
of gauge transformations. Moreover, it is in $W^{1,n}$ even across $0$,
because summing up the contributions from \re{uksmall} and \re{oksmall}
for each $k$ gives
\begin{equation*}
  \|\Omega\|^n_{W^{1,n}(B_{1/64})}
  \le 
  \sum_{k=0}^\infty\|\Omega\|^n_{W^{1,n}(Q_k)}
  \le
  C\sum_{k=0}^\infty\del_k^n
  \le 
  C\del^2,
\end{equation*}
where here, the last estimate follows from~\eqref{sum_delta}.
This means that we can interpret $\Omega$ as a $W^{1,n}$-connection form on the
trivial bundle over $B_{1/64}$, hence we have removed the singularity of
the bundle. A final application of the Uhlenbeck gauge
theorem~\ref{uhl_n} provides us with another gauge that transforms
$\Omega$ into the desired $W^{n-1,2}$-connection on the trivial bundle
over $B_{1/128}$. \qed

\section{Euler-Lagrange equations}\label{sec:euler}

The bi-Yang-Mills equation has been computed by \cite{BU}, to the effect
that the Euler-Lagrange equation of $\int_M|d_A^*F_A|^2\,dx$ is of the form
\[
  d_A^*d_Ad_A^*F_A=d_A^*F_A\#F_A,
\]
where here ``$\#$'' is used as before. We will generalize this here,
first for $\int_M|\das{n-2}F_A|^2\,dx$ and then for $Y_n$, and for $Z_n$ by
analogy. 
As a preparation, we compute

\begin{lemma}\la{1st-var}
Let $M$ be a manifold 
of dimension $m$, and $n\ge2$ be given. Let $A_t$ be a smooth
$1$-parameter family of $W^{n-1,2}$-connections on a principal bundle
$P$ over $M$ and $A:=A_0$, $\a:=\ddt A_t$. Then we have
\begin{align}
  \label{ddtdF}
  \ddt\dast{n-2}F_{A_t}
  &=\da{n-1}\a-(-1)^m\sum_{k=0}^{\lfloor\frac{n-3}2\rfloor}\da{n-3-2k}
  *[\a,*\das{2k}F_A]\nn\\
  &\quad{}+\sum_{k=0}^{\lfloor\frac{n-4}2\rfloor}\das{n-4-2k}[\a,\das{2k+1}F_A].
\end{align}
\end{lemma}

{\bf Proof.} We proceed by induction over $n$. For $n=2$, we compute
\[
  \ddt F_{A_t}=\ddt(dA_t+\tf12[A_t,A_t])=d\a+[A,\a]=d_A\a.
\]
For the inductive step, we observe that for $n\ge2$ even,
\begin{align*}
  \ddt\dast{n-1}F_{A_t}&=\ddt(d_{A_t}^*\dast{n-2}F_{A_t})\\
  &=d_A^*\Big(\ddt\dast{n-2}F_{A_t}\Big)+\Big(\ddt d_{A_t}^*\Big)\das{n-2}F_A\\
  &=d_A^*\Big(\ddt\dast{n-2}F_{A_t}\Big)-(-1)^m*[\a,*\das{n-2}F_A],
\end{align*}
while for odd $n\ge3$, similarly,
\[
  \ddt\dast{n-1}F_{A_t}=d_A\Big(\ddt\dast{n-2}F_{A_t}\Big)+[\a,\das{n-2}F_A].
\]
The assertion now follows easily.\qed

\begin{lemma}[Euler-Lagrange equations]\la{el-lem-diff}
Let $M$ be of any dimension $m$, and $n\ge2$ be given. 
The Euler-Lagrange equation for $\int_M|d^{*\wedge n-2}_AF_A|^2\dx$ is given
by
\[
  \das{2n-3}F_A=\sum_{\ell=0}^{n-3}\das{2n-5-\ell}F_A\#\das{\ell}F_A,
\]
with bilinear forms noted as above. The Euler-Lagrange equation for $Y_n$ reads
\[
  \das{2n-3}F_A+\tf{n}2\,d_A^*(|F_A|^{n-2}F_A)
  =\sum_{\ell=0}^{n-3}\das{2n-5-\ell}F_A\#\das{\ell}F_A.
\] 
\end{lemma}

{\bf Proof.} 
Assume that $A$ is a critical point of $\int_M|d^{*\wedge
  n-2}_AF_A|^2\dx$ and $A_t$ a variation of $A$ as in the Lemma \ref{1st-var}.
The lemma then yields
\eas
  0&=&\frac12\ddt\int_M|\dast{n-2}F_{A_t}|^2\,dx\\
  &=&\int_M\<\das{n-2}F_A\,,\,\da{n-1}\a\>\,dx\\
  &&{}-(-1)^m\sum_{k=0}^{\lfloor\frac{n-3}2\rfloor}\int_M\<\das{n-2}F_A\,,\,
  \da{n-3-2k}*[\a,*\das{2k}F_A]\>\,dx\\
  &&{}+\sum_{k=0}^{\lfloor\frac{n-4}2\rfloor}\int_M\<\das{n-2}F_A\,,\,
  \das{n-4-2k}[\a,\das{2k+1}F_A]\>\,dx\\
  &=&\int_M\<\das{2n-3}F_A\,,\,\a\>\,dx\\
  &&{}+\sum_{k=0}^{\lfloor\frac{n-3}2\rfloor}\int_M\<*\das{2n-5-2k}F_A\,,\,
  [\a,*\das{2k}F_A]\>\,dx\\
  &&{}+\sum_{k=0}^{\lfloor\frac{n-4}2\rfloor}\int_M\<\das{2n-6-2k}F_A\,,\,
  [\a,\das{2k+1}F_A]\>\,dx,
\eeas
from which we read off the assertion. The additional term for $Y_n$ is
straightforward.\qed

While Lemma \ref{el-lem-diff} writes the Euler-Lagrange equation in a rather
simple form, what we need is some divergence form to make sense of weak
solutions.

\begin{lemma}[Euler-Lagrange equations, divergence form]\la{el-lem}
Let $M$ be of any dimension $m$, and $n\ge2$ be given. 
Then we have smooth coefficient forms
$P_k$ ($k\in\{n+1,\ldots2n-1\}$) such that $P_k[A]$ depends on
$A,DA,\ldots,D^{n-1}A$ and is a $1$-form if $k$ is odd, a $2$-form if
$k$ is even, and that the Euler-Lagrange equation for 
$\int_M|d^{*\wedge n-2}_AF_A|^2\dx$ reads
\eq{eula1}
  \dn{(2n-2)}A=\sum_{k=n+1}^{2n-1}\dnx{(2n-1-k)}P_k[A].
\eeq
If $A$ is in Coulomb gauge ($d^*A=0$), this can be re-written as
\eq{eula2}
  (-1)^{n-1}\Delta^{n-1}A=\sum_{k=n+1}^{2n-1}\dnx{(2n-1-k)}P_k[A].
\eeq
Here we wrote $\dnx{i}$ for $d^{\wedge i}$ if operating on 1-forms,
respectively for $d^{*\wedge i}$ if acting on 2-forms,
cf. Section~\ref{sec:diff-forms} for the definition. 
Each $P_k[A]$ satisfies
\eq{eula3}
  |P_k[A]|\le C\Big(1+\sum_{j=0}^{n-1}|D^jA|^{\frac{k}{j+1}}\Big)\qquad \mbox{for
  }k\in\{n+1,\ldots,2n-1\}. 
\eeq
For the functional $Y_n$, the same formulae \re{eula1}--\re{eula3} hold, 
with different forms $P_{2n-2}$ and $P_{2n-1}$ of the same structure.

For the functional $Z_n$, \re{eula1} and \re{eula2} have to be modified
and read
\eq{eula1a}
  \Delta^{n-2}d^*dA=\sum_{k=n+1}^{2n-1}(D^*)^{2n-1-k}P_k[A]
\eeq
for every gauge, and
\eq{eula2a}
  -\Delta^{n-1}A=\sum_{k=n+1}^{2n-1}(D^*)^{2n-1-k}P_k[A]
\eeq
in Coulomb gauge. In this case, the $P_k[A]$ are sections of 
$\otimes^kT^*M\otimes\gg_P$, and again, the estimates \re{eula3} hold.
\end{lemma}

{\bf Proof.} 
Again we apply Lemma \ref{1st-var} to critical points of $\int_M|d^{*\wedge
  n-2}_AF_A|^2\dx$, but this time after differentiating through all the
terms of the lemma.
\begin{align}
  \label{eu4}
  0&=\frac12\ddt\int_M|\dast{n-2}F_{A_t}|^2\,dx\\\nn
  &=\int_M\Big\<\das{n-2}F_A\,,\, \da{n-1}\a\Big\>\dx\\
  &\quad+\int_M\Big\<\das{n-2}F_A\,,\,\sum_{j\in K_n}\da{j_1-1}\a\#
  (\das{j_2-2}F_A\#\ldots\#\das{j_\ell-2}F_A)\Big\>\,dx.\nn
\end{align}
Here
\[
  K_n:=\{j=(j_1,\ldots,j_\ell):\;\ell\ge2,\;j_1\ge1,\;
    j_2\ge\ldots\ge j_\ell\ge2,\;j_1+\ldots+j_\ell\le n\}.
\]
The only reason that $K_n$ contains sets with $j_1+\ldots+j_\ell<n$ is that
such terms appear when the forms ``$\#$'' themselves are differentiated.

We note that by the definition of $K_n$, we have
$j_1-1\in\{0,\ldots,n-3\}$, note, however, that in the first term on the
right-hand side of \re{eu4} there is an $(n-1)$-th derivative of $\a$. 
Therefore, the above equation can be transformed into 
\begin{equation}\label{weak-eula}
  \int_M \<d^{\wedge n-1}A,d^{\wedge n-1}\alpha\>
         +\sum_{k=n+1}^{2n-1}\<d^{\wedge(2n-1-k)}\alpha,P_k[A]\>\dx=0
\end{equation}
for functions $P_k[A]$ of the form
\begin{equation*}
  P_k[A]=\sum_{j_1+\ldots+j_\ell\le k}D^{j_1-1}A\#\ldots\#D^{j_\ell-1}A.
\end{equation*}
The claimed estimate \re{eula3} follows by applying Young's inequality
with exponents $\frac k{j_1},\ldots,\frac k{j_\ell}$ (and then again Young
if the exponents are too small). 
The equation~\re{weak-eula} is the weak formulation of the
assertion \re{eula1}. If $A$ is in Coulomb gauge, then
$(-1)^{n-1}\Delta^{n-1}A=d^{\wedge(2n-2)}A$, which yields \re{eula2}.

For $Y_n$, the additional term $\frac{n}2\,d_A^*(|F_A|^{n-2}F_A)$
contributes to both $P_{2n-2}$ and $P_{2n-1}$ provided $n>2$, while
for $n=2$, it coincides with the term $d_A^{*\wedge 2n-3}F_A$ that was
already treated above.

Finally, for $Z_n$, everything works quite similarly, replacing alternating
powers of $d_A$ and $d_A^*$ by suitable combinations of $D_A$ and $D_A^*$.
Note that the leading term in Lemma \ref{el-lem-diff} becomes 
$d_A^*(D_A^*)^{n-2}D_A^{n-2}F_A$.\qed

\section{Smoothness of weak solutions}\label{sec:regularity}

Next we
apply elliptic bootstrap arguments to establish smoothness of solutions
for the Euler-Lagrange equations of $Y_n$ or $Z_n$, or many similar functionals.
Note that in particular our result includes regularity of weakly
bi-Yang-Mills connections in the sense of \cite{BU}, in dimensions $m\le6$.

\begin{theorem}\label{thm:regularity}
  Assume that $A\in W^{n-1,2}(B_1^m,\wedge^1\R^m\otimes\gg)$ is a
  weak solution of \re{eula1} or \re{eula1a}, and that it is in Coulomb
  gauge, i.\,e. $d^*A=0$. We suppose that (\ref{eula3}) is in force
  and that $m\le 2n$. Then there holds $A\in
  C^\infty(B_1,\wedge^1\R^m\otimes\gg)$. 
\end{theorem}

\textbf{Proof. }The proof consists of three steps. We begin with\\\strut 

\textbf{Step 1: A Morrey-type estimate. }
The goal in this step is to establish a Morrey-type estimate of the type
\begin{equation}\label{claim_step1}
  \sum_{j=0}^{n-1}\sup_{B_\rho(y)\subset B_R}\rho^{-2\alpha}\bigg(\rho^{2n-m}\int_{B_\rho}|D^jA|^{\frac{2n}{j+1}}\dx\bigg)^{\frac{j+1}n}<\infty
\end{equation}
for any $\alpha\in(0,1)$ and every $R>0$. We note that in the
subcritical case $m<2n$, this estimate is trivially satisfied at least for some
$\alpha\in(0,1)$, but this first step of the proof is crucial in the critical case
$m=2n$. 

We temporarily fix $y\in B_R$.
For an $\eps_0\in(0,1)$ to be fixed later, we choose $r_0\in(0,1)$ 
small enough to achieve $B_{r_0}(y)\subset B_R$ and 
\begin{equation}
  \label{eps0}
  \Phi(r_0):=\sum_{j=0}^{n-1}\bigg(r_0^{2n-m}\int_{B_{r_0}(y)}|D^jA|^{\frac{2n}{j+1}}\dx
    \bigg)^{\frac{j+1}n}<\eps_0.
\end{equation}
For the remainder of
this first step, we always consider a radius $r\in(0,r_0)$ and abbreviate
$B_r:=B_r(y)\subset B_R$. 
We decompose the solution into $A=A_0+A_1$, where
$A_1$ is the $\Delta^{n-1}$-polyharmonic function (component-wise) with
the same boundary values as $A$. Since $A_0\in W^{n-1,2}_0(B_r)$,
$A$ solves the Euler-Lagrange equation and $\Delta^{n-1}(A-A_0)=0$,  there holds
\begin{align*}
  &r^{2n-m}\int_{B_r}|D^{n-1}A_0|^2\dx\\
  &\quad\le Cr^{2n-m}\sum_{k+1}^{2n-1}\int_{B_r} |P_k[A]|\,|D^{2n-1-k}A_0|\dx\\
  &\quad\le C\sum_{k+1}^{2n-1}\bigg(r^{2n-m}\int_{B_r}
  |P_k[A]|^{\frac{2n}k}\dx\bigg)^{\frac
    k{2n}}\bigg(r^{2n-m}\int_{B_r}|D^{2n-1-k}A_0|^{\frac{2n}{2n-k}}\dx\bigg)^{\frac{2n-k}{2n}}\\
  &\quad\le C\sum_{k+1}^{2n-1}\bigg(r^{2n-m}\int_{B_r}
  |P_k[A]|^{\frac{2n}k}\dx\bigg)^{\frac k{2n}} 
  \bigg(r^{2n-m}\int_{B_r}|D^{n-1}A_0|^2\dx\bigg)^{\frac12}.
\end{align*}
Here, we used the Sobolev embedding $W^{n-1,2}_0\hookrightarrow
W^{2n-1-k,\frac{2n}{2n-k}}$ for the last step, which holds in any
dimension $m\le 2n$. Re-absorbing the last
integral on the left-hand side and using the Poincar\'e-Sobolev inequality, we arrive at
\begin{align}
  \label{comparison}
  \sum_{j=0}^{n-1}\bigg(r^{2n-m}
      \int_{B_r}|D^jA_0|^{\frac{2n}{j+1}}\dx\bigg)^{\frac{j+1}n}  
  &\le Cr^{2n-m}\int_{B_r}|D^{n-1}A_0|^2\dx\\\nonumber
  &\le C\sum_{k+1}^{2n-1}
      \bigg(r^{2n-m}\int_{B_r} |P_k[A]|^{\frac{2n}k}\dx\bigg)^{\frac k{n}}.
\end{align}
Since $D^jA_1\in C^\infty(B_r)$ is polyharmonic for every $j\in\N$, we get similarly as in \cite[Lemma 6.2]{GS}
 \begin{align}\label{polyharmonic}
   \mint_{B_\rho}|D^jA_1|^{\frac{2n}{j+1}}\dx
   &\le
    C\|D^jA_1\|_{L^\infty(B_{r/2})}^{\frac{2n}{j+1}}\\\nn
   &\le 
   C\bigg(\mint_{B_r}|D^jA_1|\dx\bigg)^{\frac{2n}{j+1}}
   \le
   C\mint_{B_r}|D^jA_1|^{\frac{2n}{j+1}}\dx
 \end{align}
for all $\rho\in(0,\frac r2)$. Now we first apply \re{polyharmonic}
with $j\in\{0,\ldots,n-1\}$ and then \re{comparison}, with the result
\begin{align*}
  &\rho^{2n-m}\int_{B_\rho}|D^jA|^{\frac{2n}{j+1}}\dx\\
  &\quad\le C\rho^{2n-m}\int_{B_\rho}|D^jA_0|^{\frac{2n}{j+1}}\dx
      +C\rho^{2n}\mint_{B_\rho}|D^jA_1|^{\frac{2n}{j+1}}\dx\\\nn
  &\quad\le
  C\rho^{2n-m}\int_{B_\rho}|D^jA_0|^{\frac{2n}{j+1}}\dx
  +C\rho^{2n}\mint_{B_r}|D^jA_1|^{\frac{2n}{j+1}}\dx\\\nn
  &\quad
   \le Cr^{2n-m}
   \int_{B_r}|D^jA_0|^{\frac{2n}{j+1}}\dx
   +C\Big(\frac\rho r\Big)^{2n}r^{2n-m}\int_{B_r}|D^jA|^{\frac{2n}{j+1}}\dx\\\nn
  &\quad\le C\sum_{k+1}^{2n-1}\bigg(r^{2n-m}\int_{B_r}
  |P_k[A]|^{\frac{2n}k}\dx\bigg)^{\frac k{j+1}}
   +C\Big(\frac\rho r\Big)^{2n}r^{2n-m}\int_{B_r}|D^jA|^{\frac{2n}{j+1}}\dx
\end{align*}
and after taking roots, we infer
\begin{align}  \label{decay1}
  &\bigg(\rho^{2n-m}\int_{B_\rho}|D^jA|^{\frac{2n}{j+1}}\dx\bigg)^{\frac{j+1}n}\\\nn
  &\qquad\le C\sum_{k+1}^{2n-1}\bigg(r^{2n-m}\int_{B_r}
  |P_k[A]|^{\frac{2n}k}\dx\bigg)^{\frac kn}\\\nonumber
  &\qquad\qquad+C\Big(\frac\rho r\Big)^{2(j+1)}\bigg(r^{2n-m}\int_{B_r}|D^jA|^{\frac{2n}{j+1}}\dx\bigg)^{\frac{j+1}n}
\end{align}
for every $j\in\{0,\ldots,n-1\}$. 
By the bounds \re{eula3} for the functions $P_k$, we can estimate 
\begin{align*}
  \bigg(r^{2n-m}\int_{B_r}|P_k[A]|^{\frac{2n}k}\dx\bigg)^{\frac kn}
  &\le C \bigg(\sum_{j=0}^{n-1}r^{2n-m}
  \int_{B_r}|D^jA|^{\frac{2n}{j+1}}\dx+r^{2n}\bigg)^{\frac kn}\\
  &\le C\sum_{j=0}^{n-1}\Phi(r)^{\frac k{j+1}}+Cr^{2k}
  \le C\eps_0^{1/n}\Phi(r)+Cr^{2(n+1)}.
\end{align*}
In the last step, we used the fact that
$k\ge n+1$ and $j+1\le n$ in each term of the preceding sum and the
property \re{eps0}. Combining the two preceding estimates for
$\rho=\theta r$ and summing over $j=\{0,\ldots,n-1\}$, we arrive at
\begin{equation*}
  \Phi(\theta r)
  \le C\big(\eps_0^{1/n}+\theta^2\big)\Phi(r)+Cr^{2(n+1)}
\end{equation*}
for every $\theta\in(0,\frac12)$.
For an arbitrary $\alpha\in(0,1)$, we now choose $\theta\in(0,\frac12)$ such
that $C\theta^{2}\le \frac12\theta^{2\alpha}$ and then $\eps_0\in
(0,1)$ so small that $C\eps_0^{1/n}\le\frac12\theta^{2\alpha}$. With
these choices of parameters, the preceding bound reads
$\Phi(\theta r) \le \theta^{2\alpha}\Phi(r)+Cr^{2(n+1)}$,
which can be iterated to give
\begin{equation}\label{Morrey}
  \Phi(\rho) \le C(\alpha)\Big[\Big(\frac\rho
  {r_0}\Big)^{2\alpha}\Phi(r_0)+\rho^{2\alpha}\Big]
\end{equation}
for every $\rho\in(0,r_0)$, where $\alpha\in(0,1)$ can be chosen
arbitrarily. This proves the assertion \re{claim_step1}.\\

\textbf{Step 2: $C^{k,\alpha}$-regularity for all $k<n-1$. }
In this step, we wish to improve the Morrey estimate from the
preceding step to 
\begin{equation}
  \label{claim_step2}
  \sum_{j=[\alpha]}^{n-1}\sup_{B_\rho(y)\subset B_R}\rho^{-2\alpha}
  \bigg(\rho^{2n-m}\int_{B_\rho(y)}|D^j A|^{\frac{2n}{j+1}}\dx\bigg)^{\frac{j+1}n}<\infty
\end{equation}
for every $\alpha\in(0,n)$ and every $R<1$. Here, $[\alpha]$ denotes
the largest integer smaller than or equal to $\alpha$. If \re{claim_step2} did not hold,
we would have
\begin{equation*}
  \overline\alpha
  :=\sup\{\alpha\in(0,n)\,:\, \mbox{\re{claim_step2} is valid for }\alpha\}<n.
\end{equation*}
We note that according to
the Morrey-type estimate \re{Morrey} established in Step~1, the claim
holds true for every $\alpha\in(0,1)$, and consequently,
$\overline\alpha\ge1$. In order to derive a contradiction, we choose
an exponent $\alpha\not\in\N$ with
$0<\alpha<\overline\alpha<\alpha(1+\frac1n)<n$. Since
$\alpha<\overline\alpha$, we know that \re{claim_step2} is
valid for this value of $\alpha$, which implies in particular that
\begin{equation*}
  \mint_{B_\rho(y)}\big|D^{[\alpha]}A\big|^{\frac{2n}{[\alpha]+1}}\dx
  \le \Hat C\rho^{\alpha\frac{2n}{[\alpha]+1}-2n}  
  \qquad\mbox{for all }B_\rho(y)\subset B_R
\end{equation*}
for every $R<1$. Here and in the rest of the proof, we follow the
convention to write $C$ for constants that depend at most on $n$ and
$\gg$ and $\Hat C$ for constants that may additionally depend on
$R$ and $A$. Since $\alpha$ was chosen as non-integer, the exponent
on the right-hand side of the preceding estimate satisfies 
$
  \alpha\frac{2n}{[\alpha]+1}-2n>[\alpha]\frac{2n}{[\alpha]+1}-2n
  =-\frac{2n}{[\alpha]+1}.
$ 
Therefore,
the Dirichlet growth theorem implies $A\in C^{j,\gamma}_\loc(B_1)$ for every
$j\in\{0,\ldots,[\alpha]-1\}$ and some $\gamma>0$, which implies in particular
\begin{equation}\label{low_order}
  \sum_{j=0}^{[\alpha]-1}
  \int_{B_r}|D^j A|^{\frac{2n}{j+1}}\dx
  \le \Hat Cr^m,
\end{equation}
for every ball $B_r\subset B_R$, where the constant $\Hat C$ might depend on the $C^{[\alpha]-1}$-norm of $A$ on
$B_R$. Now we use the estimates
\re{eula3} for the functions $P_k$, together with \re{claim_step2} and
\re{low_order} in order to deduce that there holds
 for every $k\in\{n+1,\ldots,2n-1\}$
\begin{align}\label{PkEst}
  \bigg(r^{2n-m}\int_{B_r}|P_k[A]|^{\frac{2n}k}\dx\bigg)^{\frac kn}
  &\le C\bigg(\sum_{j=0}^{n-1}r^{2n-m}\int_{B_r}|D^jA|^{\frac{2n}{j+1}}\dx+r^{2n}\bigg)^{\frac kn}\\\nn
  &\le \Hat Cr^{2k}+\Hat C\sum_{j=[\alpha]}^{n-1}r^{2\alpha\frac k{j+1}}\\\nn
  &\le \Hat Cr^{2(n+1)}+\Hat Cr^{2\alpha(1+\frac1n)}
   \le \Hat Cr^{2\alpha(1+\frac1n)}
\end{align}
since $r<1$ and $\alpha(1+\frac1n)<n$ by the choice of $\alpha$. 
Combining this with the excess decay estimate \re{decay1}, we infer
for every $j\in\{0,\ldots,n-1\}$
\begin{align}  \label{decay2}
  &\bigg(\rho^{2n-m}\int_{B_\rho}|D^jA|^{\frac{2n}{j+1}}\dx\bigg)^{\frac{j+1}n}\\\nonumber
  &\qquad\le C\Big(\frac\rho 
  r\Big)^{2(j+1)}\bigg(r^{2n-m}\int_{B_r}|D^jA|^{\frac{2n}{j+1}}\dx\bigg)^{\frac{j+1}n}
  +\Hat Cr^{2\alpha(1+\frac1n)}  
\end{align}
for every $\rho\in(0,\frac r2)$. As long as
$2\alpha(1+\frac1n)<2(j+1)$, i.\,e. for $j>\alpha(1+\frac1n)-1$, we
can iterate this to get
\begin{equation*}
  \bigg(\rho^{2n-m}\int_{B_\rho(y)}|D^jA|^{\frac{2n}{j+1}}\dx\bigg)^{\frac{j+1}n}
  \le \Hat C\rho^{2\alpha(1+\frac1n)},
\end{equation*}
for every $B_\rho(y)\subset B_R$, where $\Hat C$ as before might depend on $A$ 
and on $R$, but not on
$\rho$. But this
implies \re{claim_step2} for $\alpha(1+\frac1n)>\overline\alpha$ instead of
$\alpha$, which is a contradiction to the choice of
$\overline\alpha$. We conclude that the claim \re{claim_step2} holds true
for every $\alpha\in(0,n)$. This implies in particular $A\in
C^{n-2,\alpha}_\loc(B_1)$ by the Dirichlet growth theorem. 
\\
 
\textbf{Step 3: $C^{n-1,\alpha}$-regularity and conclusion. }
In this last step, we wish to prove an excess decay estimate for the
Campanato-type excess
\begin{equation*}
  \Psi(\rho):=\rho^{2n-m}\int_{B_\rho(y)}|D^{n-1}A-(D^{n-1}A)_{y,\rho}|^2\dx,
\end{equation*}
for some $y\in B_1$ and $\rho<1-|y|$. 
To this end, we use the results from the preceding step, which imply
in particular
\begin{align*}
  \sum_{j=0}^{n-2}\int_{B_r}|D^jA|^{\frac{2n}{j+1}}\dx\le
  \Hat Cr^m
  \qquad\mbox{and}\qquad
  r^{2n-m}\int_{B_r}|D^{n-1}A|^2\dx\le \Hat C r^{2\alpha}
\end{align*}
for any $\alpha\in(n-1,n)$ and any ball $B_r=B_r(y)\subset B_R$, with constants that might depend on $A$
and $R$. In what follows, we fix a value
$\alpha\in(\frac{n^2}{n+1},n)$.  Combining the preceding estimate with the bounds \re{eula3} for $P_k$, we infer
similarly as in \re{PkEst} 
\begin{align}\label{PkEst2}\nn
  \bigg(r^{2n-m}\int_{B_r}|P_k[A]|^{\frac{2n}k}\dx\bigg)^{\frac kn}
  &\le C\bigg(r^{2n}+\sum_{j=0}^{n-1}r^{2n-m}\int_{B_r}|D^jA|^{\frac{2n}{j+1}}\dx\bigg)^{\frac kn}\\
  &\le \Hat Cr^{2k}+\Hat Cr^{2\alpha\frac kn}\le \Hat Cr^{2\alpha\frac{n+1}n}
\end{align}
for every $k\in\{n+1,\ldots,2n-1\}$. 
We again consider the decomposition $A=A_0+A_1$ into 
a $\Delta^{n-1}$-polyharmonic function $A_1$ and 
a function $A_0\in W^{n-1,2}_0(B_r)$. 
Combining \re{PkEst2} with the bound \re{comparison}, we deduce
\begin{align}
  \label{comparison2}\nonumber
  &r^{2n-m}\int_{B_r}|D^{n-1}A_0|^2\,dx\\
  &\qquad\le C\sum_{k+1}^{2n-1}
  \bigg(r^{2n-m}\int_{B_r}|P_k[A]|^{\frac{2n}k}\dx\bigg)^{\frac kn}
  \le \Hat Cr^{2\alpha\frac {n+1}n}.
\end{align}
Our next goal is an improvement of the excess decay estimate
\re{polyharmonic} for the polyharmonic function $A_1$. 
For this aim, we observe that  
with $A_1\in C^\infty(B_r)$, also  $D^nA_1$ is polyharmonic, and thus we
can estimate, following the lines of \cite[Lemma 6.2]{GS},
 \begin{align*}
   \rho^2\mint_{B_\rho}|D^nA_1|^2\dx
   &\le
    C\rho^2\|D^nA_1\|_{L^\infty(B_{r/2})}^2\\\nn
   &\le
   C\rho^2\mint_{B_r}|D^nA_1|^2\dx\\\nn
   &\le C\Big(\frac\rho r\Big)^2\mint_{B_r}|D^{n-1}A_1-(D^{n-1}A)_{B_r}|^2\dx
 \end{align*}
for all $\rho\in(0,\frac r2)$. Combining the above estimate with
Poincar\'e's inequality, we infer
\begin{equation*}
  \mint_{B_\rho}|D^{n-1}A_1-(D^{n-1}A_1)_{B_\rho}|^2\dx
  \le C\Big(\frac\rho r\Big)^2\mint_{B_r}|D^{n-1}A_1-(D^{n-1}A_1)_{B_r}|^2\dx.
\end{equation*}
Next, we transfer this decay estimate to $A$ by means of
\re{comparison2} as follows.
\begin{align*}
  \Psi(\rho)&\le\rho^{2n-m}\int_{B_\rho}|D^{n-1}A-(D^{n-1}A_1)_{B_\rho}|^2\dx\\
  &\le 2\rho^{2n-m}\int_{B_\rho}|D^{n-1}A_1-(D^{n-1}A_1)_{B_\rho}|^2\dx
    +2\rho^{2n-m}\int_{B_\rho}|D^{n-1}A_0|^2\dx\\
  &\le C\Big(\frac\rho
  r\Big)^{2n+2}r^{2n-m}\int_{B_r}|D^{n-1}A_1-(D^{n-1}A_1)_{B_r}|^2\dx\\
  &\qquad  +2\rho^{2n-m}\int_{B_\rho}|D^{n-1}A_0|^2\dx\\
  &\le C\Big(\frac\rho r\Big)^{2n+2}\Psi(r)
      +Cr^{2n-m}\int_{B_r}|D^{n-1}A_0|^2\dx\\
  &\le C\Big(\frac\rho r\Big)^{2n+2}\Psi(r)+\Hat Cr^{2\alpha\frac {n+1}n}.
\end{align*}
Since by our choice of $\alpha$, there holds $2\alpha\frac {n+1}n\in(2n,2n+2)$, 
we can iterate the above estimate to get
\begin{equation*}
  \Psi(\rho)\le \Hat C\rho^{2\alpha\frac {n+1}n}= \Hat C\rho^{2n+2\gamma}
\end{equation*}
for some $\gamma>0$ and every $0<\rho\le r\le R-|y|$. But this implies $A\in C^{n-1,\gamma}_{\rm
  loc}(B_1)$ by Campanato's integral characterization of H\"older
continuous functions. 

Having arrived at this stage, the claim $A\in C^\infty(B_1)$ follows
from classical Schauder theory, see e.\,g. \cite[Thm. 2']{DN}, which
concludes the proof of the theorem.
\qed

\section{Existence of minimizers}\label{sec:existence}

In this part of the paper, we formulate our existence theorems for both
minimizers of $Y_n$ and $Z_n$. In the proofs, however, we mention only $Y_n$.
The modifications for $Z_n$ are straightforward.

\subsection{The critical dimension}
\label{sec:existence-critical}

In the critical dimension, we encounter the problem of a possible
bubbling phenomenon during the minimizing procedure that might result
in a change of the underlying bundle.
However, certain topological invariants are preserved.
Our considerations involve prescribing certain Chern classes of the
bundle. The relations between Chern Classes and Sobolev maps with
finite $L^n$-norm of the curvature have been explored by Uhlenbeck
\cite{Uh3}, and for a minimizing procedure similar to ours by Sedlacek
\cite{Se}.

In order to
demonstrate how certain invariants are preserved, 
we restrict ourselves to the case of
principal $SU(k)$-bundles $P$ over $M$. We recall that in this case, the Chern
classes of the associated vector bundle $P_{\C^k}:=P\times_\rho \C^k$  is given by
\begin{equation*}
  c_j(P_{\C^k})=\big[p_j(\tfrac \i{2\pi}F_A)\big]\in \HdR^{2j}(M),
\end{equation*}
Here, $p_j$ denotes the $j$-th elementary symmetric polynomial of the
eigenvalues and
$F_A$ is the curvature of any connection on $P$. Moreover, 
with $P\times_\rho \C^k$ we abbreviate the complex vector bundle associated to
$P$ by the representation $\rho:SU(k)\to GL(k,\C)$ of 
$SU(k)$ on $\C^k$. 

In order to compare the Chern classes of bundles over
$M\setminus\Sigma$, the following lemma is crucial.

\begin{lemma}\label{lem:iso_inclusion}
  Assume that $M$ is a compact manifold of dimension $m\ge 4$ and
  $\Sigma\subset M$ is a finite set. Then, the inclusion
  $\iota_0:M\setminus\Sigma\to M$ induces an isomorphism
  \begin{equation*}
    \iota_0^*:\HdR^\ell(M)\to \HdR^\ell(M\setminus\Sigma)
  \end{equation*}
  for every $\ell\in\{2,\ldots,m-2\}$. Similarly, if $B$ is the union
  of finitely many, pairwise disjoint closed balls, then the inclusion 
  $\iota_1:M\setminus B\to M$ induces an isomorphism 
  \begin{equation*}
    \iota_1^*:\HdR^\ell(M)\to \HdR^\ell(M\setminus B).
  \end{equation*}
\end{lemma}

\textbf{Proof. }
For any finite set $\Sigma\subset M$ we may choose a union of 
finitely many, pairwise disjoint closed balls $B\supset\Sigma$. 
Since $M\setminus B$ is a deformation retract of $M\setminus\Sigma$,
it suffices to prove the second statement of the lemma. 
To this end, we choose a union of finitely many, pairwise disjoint \emph{open} balls $\widehat
B\supset B$. The Mayer-Vietoris sequence for the open sets
$M\setminus B$ and $\widehat B$ reads
\begin{equation*}
  \ldots\to \HdR^{\ell-1}(\widehat B\setminus B)
        \to \HdR^\ell(M)
        \overset{(\iota_1^*,\iota_2^*)}\longrightarrow
        \HdR^\ell(M\setminus B)\oplus \HdR^\ell(\widehat B)
        \to \HdR^\ell(\widehat B\setminus B)
        \to\ldots
\end{equation*}
Here, $\widehat B\setminus B$ is homotopy equivalent to $N$ spheres of
dimension $m-1>\ell$ and $\widehat B$ is the union of $N$ balls of dimension
$m>\ell$. Therefore, we have
\begin{equation*}
  \HdR^\ell(\widehat B)=0=\HdR^{\ell-1}(\widehat B\setminus B)=\HdR^\ell(\widehat B\setminus B)
\end{equation*}
for all $\ell\in\{2,\ldots,m-2\}$. Plugging this into
the Mayer-Vietoris sequence stated above, we infer that
$\iota_1^*:\HdR^\ell(M)\to\HdR^\ell(M\setminus B)$ is an isomorphism.
\qed

The following lemma will ensure that the Chern classes are preserved
under weak $L^p$-convergence. 
The proof is a slight modification of the arguments in \cite[Cor. 5.2]{ISS}.  

\begin{lemma}\label{lem:closed_cohom_class}
  Let $M$ be a regular open subset of a smooth compact manifold, $k\in\N$ and $p>1$.
  We consider a $k$-form $\phi_0\in C^\infty(\overline
  M,\wedge^kT^*M)$ and a sequence
  $\phi_i\in[\phi_0]\cap L^p(M,\wedge^kT^*M)$ with $\phi_i\wto\phi\in
  C^\infty(M,\wedge^kT^*M)$ weakly in $L^p$. Then $\phi\in[\phi_0]$.
 \end{lemma}

\textbf{Proof. }
  Since $\phi_i\in[\phi_0]$, there are smooth $(k-1)$-forms $\omega_i$
  with 
  \begin{equation*}
    d\omega_i=\phi_i-\phi_0\wto \phi-\phi_0
    \qquad\mbox{weakly in $L^p$, as $i\to\infty$}. 
  \end{equation*}
  By an approximation argument, we may assume $\omega_i\in
  C^\infty(\overline M,\wedge^{k-1}T^*M)$ for all $i\in\N$. 
  Following the strategy in \cite[Cor. 5.2]{ISS}, we can find a
  $(k-1)$-form $\omega\in W^{1,p}(M,\wedge^{k-1}T^*M)$ with 
  $d\omega_i\wto d\omega$ weakly in $L^p$ as $i\to\infty$, and in particular 
  \begin{equation*}
    d\omega=\phi-\phi_0.
  \end{equation*}
  Moreover, in \cite[Cor. 5.2]{ISS} the $(k-1)$-form $\omega$ is constructed as the weak limit
  of coexact forms in $W^{1,p}(M,\wedge^{k-1}T^*M)$, which implies
  $d^*\omega=0$. We deduce
  \begin{equation*}
    \Delta\omega=d^*d\omega=d^*(\phi-\phi_0)\in C^\infty(M,\wedge^{k-1}T^*M),
  \end{equation*}
  and elliptic regularity theory yields $\omega\in
  C^\infty(M,\wedge^{k-1}T^*M)$. We thereby arrive at  
  $\phi=\phi_0+d\omega\in[\phi_0]$, 
  which completes the proof of the lemma.
\qed

For the formulation of the minimization problem that we wish to solve, 
we fix a reference bundle $P_0$ over a manifold $M$ 
of dimension $m=2n$.
We want to prescribe the Chern
classes $c_j^0:=c_j((P_0)_{\C^k})\in H^{2j}(M)$ for $j\in
\{1,\ldots,n-1\}$. More precisely, we define the class of admissible
bundles by
\begin{equation*}
 \mathcal{P}_*(M,P_0):=\{P\in \mathcal{P}(M)\,:\,
                   c_j(P_{\C^k})=c_j^0 \mbox{ for }j=1,\ldots,n-1\},
\end{equation*}
where $\mathcal{P}(M)$ denotes the class of smooth
principal $SU(k)$-bundles over $M$. In order to emphasize the
occurrence of different bundles during the minimization process, we
will frequently use the notation
\begin{equation*}
  Y_n(A,P):= Y_n(A)
\end{equation*}
if $A$ is a smooth connection on the bundle $P\in\mathcal{P}_*(M,P_0)$. \\

Having introduced the setting, we are ready to state our existence
result in the critical case. 
\begin{theorem}\label{thm:existence}
  Let $M$ be a compact manifold of dimension $m=2n$, and $A_0$ a
  smooth reference connection on a principal $SU(k)$-bundle $P_0$ over
  $M$ as above. 
  Then there is a principal
  $SU(k)$-bundle $P\in \mathcal{P}_*(M,P_0)$ and a smooth
  connection $A$
  on $P$ that minimizes the functional $Y_n$ (or $Z_n$) in the class of smooth
  connections on bundles in $\mathcal{P}_*(M,P_0)$.
\end{theorem}

\textbf{Proof. }
The proof is divided into several steps. We start with\\

\textbf{Step 1: Convergence in Uhlenbeck gauges. }
We choose a minimizing sequence of smooth connections $A_i$ on
bundles $P_i\in\mathcal{P}_*(M,P_0)$ for the functional $Y_n$. From
Theorem~\ref{WA-estimates} we deduce
\begin{equation}\label{sup-WA-estimate}
  \sup_{i\in\N}\,\|D_{A_i}^\ell F_{A_i}\|_{L^{2n/(\ell+2)}(M)}<\infty
  \qquad\mbox{for }\ell=0,\ldots,n-2.
\end{equation}
Writing $\mu_M$ for the Riemannian measure on $M$,  
we define a sequence
of Radon measures on $M$ by
\begin{equation*}
  \mu_i:=\mu_M\edge\big(|\dasi{n-2}F_{A_i}|^2+|F_{A_i}|^n\big).
\end{equation*}
The minimizing property of $A_i$ implies
\begin{equation*}
  \sup_{i\in\N}\mu_i(M)=\sup_{i\in\N} Y_n(A_i,P_i)<\infty.
\end{equation*}
Therefore, we can find a Radon measure $\mu$ on $M$ such that
$\mu_i\wsto\mu$ weakly* in the space of Radon measures, possibly after
extracting a subsequence. With the
constant $\kappa>0$ from Theorem~\ref{uhl_n}, we define the
singular set of the limit bundle by
\begin{equation*}
  \Sigma:=\{x\in M\,:\, \mu(\{x\})\ge\kappa\}
\end{equation*}
and let $N:=\#\Sigma$.
From the definition of $\Sigma$ and $\mu$, it is evident that for every $x\in
M\setminus\Sigma$, we can find a ball $U\ni x$ such that
\begin{equation}\label{choice-U}
  \limsup_{i\to\infty}\mu_i(2U)\le \mu(\overline{2U})<\kappa.
\end{equation}
Here and in what follows, we use the notation $2U$ for the ball with
double radius and the same center as $U$.
We choose a countable cover $U_\alpha$, $\alpha\in\N$, of balls with the above
property and $\cup_{\alpha\in\N}U_\alpha =M\setminus\Sigma$.
In particular, the property~(\ref{choice-U}) implies that for every $\alpha\in\N$ and every
sufficiently large $i\ge i_0(\alpha)$, we have
\begin{equation}\label{choice-Ualpha}
  \int_{2U_\alpha}|\dasi{n-2}F_{A_i}|^2+|F_{A_i}|^n\dx<\kappa.
\end{equation}
Therefore, the Gauge Theorem~\ref{uhl_n}
provides us with gauge transformations
$u_{i,\alpha}\in W^{n,2}(U_\alpha,G)$ such that
$A_i^{\alpha}:=u_{i,\alpha}^*A_i$ satisfies
$d^*A_i^{\alpha}=0$ on $U_\alpha$ and moreover
\begin{align*}
  &\sup_{i\ge i_o(\alpha)} \sum_{\ell=0}^{n-1}\|D^\ell A_i^{\alpha}\|_{L^{2n/(\ell+1)}(U_\alpha)}\\
  &\qquad\le
  C\sup_{i\ge i_o(\alpha)}\Big(\|D_{A_i}^{n-2}F_{A_i}\|_{L^2(2U_\alpha)}+\|F_{A_i}\|_{L^n(2U_\alpha)}\Big)<\infty
\end{align*}
for every $\alpha\in\N$, where the finiteness of the right-hand side
is a consequence of~(\ref{sup-WA-estimate}). 
By extraction of a subsequence of $i$
(possibly depending on $\alpha$), we can thus achieve the
convergence to a local limit connection 
$A^{\alpha}\in W^{n-1,2}(U_\alpha,\wedge^1\R^m\otimes\suk )$ 
in the sense 
\begin{equation}\label{A-converge}
  \left\{
  \begin{array}{c l}
  A_i^{\alpha}\wto A^{\alpha}&\mbox{ weakly in
  }W^{n-1,2}(U_\alpha),
  \\[1ex]
 D^\ell A_i^{\alpha}\to D^\ell A^{\alpha}&
 \begin{array}[t]{l}
  \mbox{strongly in }L^p(U_\alpha)
   \ \mbox{ and a.e.}\\
  \mbox{for all }\ell\le n-2\mbox{ and }p<\frac{2n}{\ell+1}.
\end{array}
\end{array}
\right.
\end{equation}
In particular, the limit connection is still in Uhlenbeck gauge,
i.\,e. we have
\begin{equation}\label{Uhlenbeck-limit}
  d^*A^{\alpha}=0\qquad\mbox{on }U_\alpha.
\end{equation}
For the curvature of the connection, the above convergence implies
\begin{equation}\label{F-converge}
  \left\{
  \begin{array}{c l}
  d^{*\wedge n-2}_{A_i^\alpha}F_{A_i^\alpha}\wto d^{*\wedge n-2}_{A^{\alpha}}F_{A^{\alpha}}
  &\mbox{ weakly in }L^2(U_\alpha),
    \\[1ex]
  D_{A_i^\alpha}^\ell F_{A_i^\alpha}\wto D_{A^{\alpha}}^\ell F_{A^{\alpha}}
  &\mbox{ weakly in
  }L^{2n/(\ell+2)}(U_\alpha)
  \quad\forall \ell\le n-2,
    \\[1ex]
    F_{A_i^\alpha}\to F_{A^{\alpha}}  &
  \begin{array}[t]{l}
\mbox{strongly in }L^p(U_\alpha)
\quad\forall p<n\\
  \mbox{and a.e., provided $n>2$,}
  \end{array}
\end{array} \right.
\end{equation}
as $i\to\infty$. 
The gauge transformations $u_{i,\alpha}$ define transition functions
$\phi^i_{\alpha\beta}\in W^{n,2}(U_\alpha\cap U_\beta,SU(k))$ by the identity
\begin{equation*}
  u_{i,\beta}=u_{i,\alpha}\phi^i_{\alpha\beta}\qquad\mbox{on }U_\alpha\cap U_\beta,
\end{equation*}
for all $i\in\N$ and all $\alpha,\beta\in\N$ for which $U_\alpha\cap
U_\beta\neq\varnothing$. From the transformation rule for connections,
we have 
\begin{equation}\label{transition}
  d\phi^i_{\alpha\beta}
  =
  \phi^i_{\alpha\beta}A_i^{\beta}-A_i^{\alpha}\phi^i_{\alpha\beta}
  \qquad\mbox{on }U_\alpha\cap U_\beta.
\end{equation}
Using the fact that $A_i^{\alpha}$ and $A_i^{\beta}$ are both
bounded sequences in $W^{n-1,2}$ and $SU(k)$ is compact, we
inductively deduce from~(\ref{transition}) that 
\begin{equation*}
  \sup_{i\in\N}\|D^\ell\phi^i_{\alpha\beta}\|_{L^{2n/\ell+1}}<\infty
\end{equation*}
for all $\ell\in\{0,\ldots,n\}$. 
Therefore, we can achieve the convergence to maps
$\phi_{\alpha\beta}\in W^{n,2}(U_\alpha\cap U_\beta,SU(k))$ in the
sense 
\begin{equation}\label{phi-converge}
   \left\{
  \begin{array}{c l}
  \phi^i_{\alpha\beta}\wto \phi_{\alpha\beta}&\mbox{weakly in
  }W^{n,2}(U_\alpha\cap U_\beta,SU(k))\\[1ex]
  D^\ell \phi^i_{\alpha\beta}\to D^\ell
  \phi_{\alpha\beta}&\mbox{
              almost everywhere }\forall \ell\le n-1,
\end{array}
\right.
\end{equation}
as $i\to\infty$. 
From the almost everywhere 
convergence of the $\phi^i_{\alpha\beta}$ and $A_i^\alpha$
by~(\ref{A-converge}) and (\ref{phi-converge}), we
infer that we can pass to the limit in~(\ref{transition}), with the
result
\begin{equation}
  \label{transition-limit}
  d\phi_{\alpha\beta}
  =
  \phi_{\alpha\beta}A^{\beta}-A^{\alpha}\phi_{\alpha\beta}
  \qquad\mbox{on }U_\alpha\cap U_\beta.
\end{equation}
Moreover, the almost everywhere convergence guarantees that 
the cocycle conditions
\begin{equation}\label{cocycle}
  \phi_{\alpha\beta}\phi_{\beta\gamma}=\phi_{\alpha\gamma}
  \qquad\mbox{on }U_\alpha\cap U_\beta\cap U_\gamma
\end{equation}
are preserved in the limit. We note that if the transition functions
$\{\phi_{\alpha\beta}\}_{\alpha,\beta\in\N}$ 
were of class $C^\infty$, they would define a new
principal bundle  $P$ over
$\cup_{\alpha\in\N}U_\alpha=M\setminus\Sigma$. Therefore we turn our
attention to \\

\textbf{Step 2: Regularity of the limit configuration. }
The smoothness of the local limit connections $A^{\alpha}$ will
follow from Theorem~\ref{thm:regularity} once we have established that
$A^{\alpha}$ weakly solves the Euler-Lagrange
equations~(\ref{eula1}) on $U_\alpha$. For this it suffices to prove
\begin{equation}\label{local_minimizer}
  Y_n(A^\alpha)\le Y_n(A^\alpha+\p)
\end{equation}
for all $\p\in 
C^\infty_{\rm  cpt}(U_\alpha,\wedge^1\R^m\otimes\suk)$. If this was
not the case, 
we could find $\p\in C^\infty_{\rm cpt}(U_\alpha,\wedge^1\R^m\otimes\suk)$ for which the
inequality~(\ref{local_minimizer}) 
does not hold. The main step to reach a contradiction is to prove the

\textbf{Claim. }As $i\to\infty$, we have the convergence
\begin{equation}\label{Y-Y-converge}
  Y_n(A^\alpha_i+\p)-Y_n(A^\alpha_i)
  \longrightarrow Y_n(A^\alpha+\p)-Y_n(A^\alpha)<0.
\end{equation}
For the \textbf{proof of the claim, }we begin by calculating
\begin{equation*}
  F_{A_i^\alpha+\p}=F_{A_i^\alpha}+[A_i^\alpha,\p]+d\p+\tfrac12[\p,\p],
\end{equation*}
and the same holds for $A^\alpha$ instead of $A^\alpha_i$. As a
consequence of the strong convergence $A_i^\alpha\to A^\alpha$
in $L^{n}$ according to~(\ref{A-converge}) and $\p\in
C^\infty_{\rm cpt}$, we deduce  
\begin{equation}\label{F-diff-converge}
   F_{A_i^\alpha+\p}-F_{A_i^\alpha}
   \longrightarrow
   F_{A^\alpha+\p}-F_{A^\alpha}
   \qquad\mbox{strongly in $L^n$, as $i\to\infty$.}
\end{equation}
Abbreviating
\begin{equation*}
  I(\omega,\psi):=\int_0^1n|\omega+t(\psi-\omega)|^{n-2}\,(\omega+t(\psi-\omega))\,dt
\end{equation*}
for $\omega,\psi\in\wedge^2\R^m\otimes\suk$ 
we have $|\omega|^n-|\psi|^n
=I(\omega,\psi)\cdot(\omega-\psi)$. 
As $i\to\infty$, we moreover have
\begin{equation*}
  I(F_{A_i^\alpha+\p},F_{A_i^\alpha})
  \wto
  I(F_{A^\alpha+\p},F_{A^\alpha})
  \qquad\mbox{weakly in }L^{\frac n{n-1}}(U_\alpha,\wedge^2\otimes\suk).
\end{equation*}
In fact, in the case $n>2$, this convergence holds almost everywhere
by~(\ref{F-converge}) and the weak convergence follows since 
$I(F_{A_i^\alpha+\p},F_{A_i^\alpha})$ is bounded in $L^{n/(n-1)}$ because
of $|I(F_{A_i^\alpha+\p},F_{A_i^\alpha})|\le
n(|F_{A_i^\alpha+\p}|+|F_{A_i^\alpha}|)^{n-1}$. For $n=2$, however,
the term $I(F_{A_i^\alpha+\p},F_{A_i^\alpha})$ depends linearly on the
curvature and therefore, the claim follows from the weak convergence
in~(\ref{F-converge}). 
Combining the above results, we infer
\begin{align}\label{F-F-converge}\nn
  \int_{U_\alpha}|F_{A_i^\alpha+\p}|^n-|F_{A_i^\alpha}|^n\dx
  &= \int_{U_\alpha}I(F_{A_i^\alpha+\p},F_{A_i^\alpha})\cdot
    (F_{A_i^\alpha+\p}-F_{A_i^\alpha})\dx\\\nn
  &\underset{i\to\infty}\longrightarrow
    \int_{U_\alpha}I(F_{A^\alpha+\p},F_{A^\alpha})\cdot
    (F_{A^\alpha+\p}-F_{A^\alpha})\dx\\
  &=\int_{U_\alpha}|F_{A^\alpha+\p}|^n-|F_{A^\alpha}|^n\dx.
\end{align}
Similarly, we proceed for the derivatives of the
curvature. Inductively, we calculate for every
$\ell\in\{1,\ldots,n-2\}$
\begin{align}\label{DF-DF}
    &d^{*\wedge\ell}_{A_i^\alpha+\p}F_{A_i^\alpha+\p}
    -d^{*\wedge\ell}_{A_i^\alpha}F_{A_i^\alpha}\\\nonumber
    &\qquad=\sum_{j}\sum_{k=1}^K D^{j_1-1}\p\#\ldots D^{j_k-1}\p\#
     D^{j_{k+1}-1}A_i^\alpha\#\ldots D^{j_K-1}A_i^\alpha, 
\end{align}
where the sum is taken over all tuples $j=(j_1,\ldots,j_K)$ of
$K\in\N$ naturals with $j_1+\ldots+j_K\le\ell+2$. We
note that at least one of the factors in each term of 
the above sum contains $\p$
and therefore is bounded. Moreover, we note that for every partition
$j$ occurring in the above sum for $\ell=n-2$, we have 
$\frac{j_{k+1}}{2n}+\ldots+\frac{j_K}{2n}<\frac12$. Because of the strong
convergence~(\ref{A-converge})$_2$, we thereby infer
\begin{equation}\label{DF-diff-converge}
    d^{*\wedge n-2}_{A_i^\alpha+\p}F_{A_i^\alpha+\p}
    -d^{*\wedge n-2}_{A_i^\alpha}F_{A_i^\alpha}
      \longrightarrow
   d^{*\wedge n-2}_{A^\alpha+\p}F_{A^\alpha+\p}-d^{*\wedge n-2}_{A^\alpha}F_{A^\alpha}
   \qquad\mbox{ in $L^2$,}
\end{equation}
as $i\to\infty$. Joining this with the weak convergence in~(\ref{F-converge}), we deduce
\begin{align}\label{DF-DF-converge}\nn
  &\int_{U_\alpha}|d^{*\wedge n-2}_{A_i^\alpha+\p}F_{A_i^\alpha+\p}|^2
  -|d^{*\wedge n-2}_{A_i^\alpha}F_{A_i^\alpha}|^2\dx\\\nn
  &\qquad=
  \int_{U_\alpha}(d^{*\wedge
    n-2}_{A_i^\alpha+\p}F_{A_i^\alpha+\p}+d^{*\wedge n-2}_{A_i^\alpha}F_{A_i^\alpha})\cdot
    (d^{*\wedge n-2}_{A_i^\alpha+\p}F_{A_i^\alpha+\p}-d^{*\wedge n-2}_{A_i^\alpha}F_{A_i^\alpha})\dx\\\nn
  &\qquad\underset{i\to\infty}\longrightarrow
    \int_{U_\alpha}(d^{*\wedge
      n-2}_{A^\alpha+\p}F_{A^\alpha+\p}+d^{*\wedge n-2}_{A^\alpha}
    F_{A^\alpha})
    \cdot (d^{*\wedge n-2}_{A^\alpha+\p}F_{A^\alpha+\p}
  -d^{*\wedge n-2}_{A^\alpha}F_{A^\alpha})\dx\\
  &\qquad=\int_{U_\alpha}|d^{*\wedge n-2}_{A^\alpha+\p}F_{A^\alpha+\p}|^2
    -|d^{*\wedge n-2}_{A^\alpha}F_{A^\alpha}|^2\dx.
\end{align}
We can combine the convergences~(\ref{F-F-converge})
and~(\ref{DF-DF-converge}) to yield the claim~(\ref{Y-Y-converge}). 

Now that the claim is proven, we resume Step 2.
We define new connections $\widetilde A_i$ on the bundles
$P_i$ in such a way that they agree with $A_i$ on $P_i|_{M\setminus
  U_\alpha}$ and such that
$u_{i,\alpha}^*\widetilde A_i=A_i^{\alpha}+\varphi$ on $P_i|_{U_\alpha}$. In view
of~(\ref{Y-Y-converge}), we have 
\begin{align*}
  Y_n(\widetilde A_i,P_i)-Y_n(A_i,P_i)
  &=Y_n(A^\alpha_i+\p)-Y_n(A^\alpha_i)\\
  &\underset{i\to\infty}{\longrightarrow}Y_n(A^\alpha+\p)-Y_n(A^\alpha)
  <0
\end{align*}
by the choice of $\p$ as a counterexample
to~(\ref{local_minimizer}). Therefore, in both
cases we infer 
\begin{equation*}
  \lim_{i\to\infty}Y_n(\widetilde A_i,P_i)<
  \lim_{i\to\infty}Y_n(A_i,P_i)=\inf_{(A,P)}Y_n(A,P),
\end{equation*}
where the infimum is taken over all $P\in\mathcal{P}_*(M,P_0)$ and all
smooth connections $A$ on $P$. We have thus reached the
contradiction that proves~(\ref{local_minimizer}). Lemma~\ref{el-lem}
now yields that $A^\alpha$ is a weak solution of (\ref{eula1}), and
since it is moreover in Uhlenbeck gauge according
to~(\ref{Uhlenbeck-limit}), we infer from Theorem~\ref{thm:regularity}
that $A^\alpha\in
C^\infty(U_\alpha,\wedge^1\R^m\otimes\suk)$. Furthermore, from the
transition formula~(\ref{transition-limit}) we inductively deduce
$|D^\ell\phi_{\alpha\beta}|\in L^\infty_\loc(U_\alpha\cap U_\beta)$ for all
$\ell\in\N$ and thereby $\phi_{\alpha\beta}\in C^\infty(U_\alpha\cap
U_\beta,SU(k))$. Keeping in mind the cocycle
conditions~(\ref{cocycle}), we can construct a new principal bundle
$P$ over
$\bigcup_{\alpha\in\N}U_\alpha=M\setminus\Sigma$ with the transition
functions $\{\phi_{\alpha\beta}\}$. From~(\ref{transition-limit}), 
we moreover know
\begin{equation*}
  A^\beta=\phi_{\alpha\beta}^{-1}A^\alpha\phi_{\alpha\beta}+\phi_{\alpha\beta}^{-1}d\phi_{\alpha\beta}
  \qquad\mbox{on }U_\alpha\cap U_\beta
\end{equation*}
for all $\alpha,\beta\in\N$ for which the latter set is not
empty. This means that the collection $\{A^\alpha\}$ defines a
connection $A$ on the new bundle $P$ over $M\setminus\Sigma$.\\

\textbf{Step 3: Removability of the point singularities}

We write $\Sigma=\{x_1,\ldots,x_N\}$ for the singular set of the bundle
$P$ and choose pairwise disjoint open balls $V_\ell\subset M$ 
with centers in the points $x_\ell\in\Sigma$ for $\ell\in\N$. 
Then, we can find a
finite subset $I\subset\N$ with
\begin{equation}\label{open-cover}
  M=\bigcup_{\ell=1}^N V_\ell\ \cup\ \bigcup_{\alpha\in I} U_\alpha.
\end{equation}
As a consequence of~(\ref{sup-WA-estimate}), the lower
semicontinuity of the norm with respect to weak convergence and the
gauge invariance of $D_A^\ell F_A$, we have
\begin{equation*}
  \sum_{\ell=0}^{n-2}\|D_{A}^\ell F_{A}\|_{L^{2n/(\ell+2)}(M)}<\infty.
\end{equation*}
We may thus
apply the Removable Singularity Theorem~\ref{RemSing} on each
ball $V_\ell\setminus\{x_\ell\}$. This theorem provides us with 
local $W^{n,2}_\loc$-trivializations
$v_\ell:\pi^{-1}(V_\ell)\to (V_\ell\setminus\{x_\ell\})\times SU(k)$
of the bundle $\pi:P\to M\setminus\Sigma$ and $\widehat A_\ell\in
W^{n-1,2}(V_\ell,\wedge^1\R^m\otimes\suk)$ such that 
$$
  (v_\ell)_*A=\widehat A_\ell|_{V_\ell\setminus\{x_\ell\}}
  \qquad\mbox{for }\ell=1,\ldots,N.
$$
These trivializations give rise to the transition functions
$\psi_{\alpha\ell}\in W^{n,2}(V_\ell\cap U_\alpha,SU(k))$ defined by
\begin{equation*}
  v_\ell=u_\alpha\psi_{\alpha\ell}\qquad\mbox{on }V_\ell\cap U_\alpha,
\end{equation*}
for any $\ell\in\{1,\ldots,N\}$ and $\alpha\in I$ for which
$V_\ell\cap U_\alpha\neq\varnothing$. 
By diminishing the balls $V_\ell$ is necessary, we can achieve that
the smallness assumption from the Gauge Theorem~\ref{uhl_n} is
satisfied on $2V_\ell$. 
Therefore, we may assume that the trivializations $v_\ell$
have been chosen in such a way that $d^*\widehat A_\ell=0$ holds for
every $\ell\in\{1,\ldots,N\}$. The arguments from Step~2 of the proof
therefore imply that $\widehat A_\ell$ weakly solves an Euler-Lagrange
equation of the form~\eqref{eula2} on
$V_\ell\setminus\{x_\ell\}$. A standard capacity argument then shows
that $\widehat A_\ell$ actually solves~\eqref{eula2} on the whole ball
$V_\ell$. We are thus in a position to apply the Regularity
Theorem~\ref{thm:regularity} to deduce that 
$\widehat A_\ell\in C^\infty(V_\ell,\wedge^1\R^m\otimes\suk)$. In the same manner as in
Step~2 we then infer that the transition functions $\psi_{\alpha\ell}$
are smooth. This implies that the collection of transition functions
$\{\phi_{\alpha\beta}\}\cup\{\psi_{\alpha\ell}\}$ 
relative to the open cover~(\ref{open-cover})
defines a smooth bundle $\widehat P$ over $M$, and the collection
$\{A_\alpha\}\cup\{\widehat A_\ell\}$ of local connection forms defines a
smooth connection $\widehat A$ on $\widehat P$.\\

\textbf{Step 4: The Chern classes are preserved. }
Next we wish to prove that $\widehat P$ possesses the same Chern classes as the
reference bundle $P_0$ and therefore, $\widehat P\in\mathcal{P}_*(M,P_0)$
is an admissible comparison bundle.
Since the inclusion $\iota:M\setminus\Sigma\to M$ induces an
isomorphism
$\iota^*:\HdR^{2j}(M)\to\HdR^{2j}(M\setminus\Sigma)$ 
by Lemma~\ref{lem:iso_inclusion},
it suffices to show that
\begin{equation}\label{Claim}
  c_j(P_{\C^k})=\iota^*c_j^0\in \HdR^{2j}(M\setminus\Sigma)
  \qquad\mbox{for }j=1,\ldots,n-1. 
\end{equation}
We choose a union of $N$ pairwise disjoint closed balls with centers
in the points of $\Sigma$. Since $M\setminus B$ is a deformation
retract of $M\setminus\Sigma$, it suffices to show
\begin{equation}\label{Chern-claim}
  p_j\big(\tfrac{\mathbf{i}}{2\pi}F_A\big|_{M\setminus B}\big)\in \iota_B^*c_j^0\in
  \HdR^{2j}(M\setminus B),
\end{equation}
where $\iota_B:M\setminus B\hookrightarrow M$ is the inclusion. 
Because $\overline{M\setminus B}$ is compact, it is covered by
finitely many of the open sets $U_\alpha$. 
Since the $A_i$ are smooth connections on principal bundles $P_i\in
\mathcal{P}_\ast(M,P_0)$, of which we prescribed the corresponding Chern
classes, we know
\begin{equation*}
  p_j\big(\tfrac{\mathbf{i}}{2\pi}F_{A_i}\big)
  \in 
  c_j^0\in \HdR^{2j}(M),
\end{equation*}
for all $i\in\N$ and $j=1,\ldots,n-1$. Now we
apply Lemma~\ref{lem:iso_inclusion}, which states that 
$\iota_B$ induces an isomorphism on the
cohomology groups of order $2j$ for $j=1,\ldots,n-1$. The preceding
formula thereby implies
\begin{equation}\label{Chern-i}
  p_j\big(\tfrac{\mathbf{i}}{2\pi}F_{A_i}\big|_{M\setminus B}\big)\in \iota_B^*c_j^0\in
  \HdR^{2j}(M\setminus B).
\end{equation}
In order to prove~(\ref{Chern-claim}), it therefore suffices to check
that this property is preserved in the limit $i\to\infty$. We first
consider the case $n>2$, in which the strong convergence
in~(\ref{F-converge}) implies strong convergence $F_{A_i^\alpha}\to
F_{A^\alpha}$ of the local representations of
$F_{A_i}$ in $L^p$-norm for every $p<n$. Keeping in mind that 
$p_j$ is a polynomial of order $j$ and using gauge invariance, we deduce
\begin{equation*}
  p_j\big(\tfrac{\mathbf{i}}{2\pi}F_{A_i}\big|_{M\setminus B}\big)
  \to
  p_j\big(\tfrac{\mathbf{i}}{2\pi}F_{A}\big|_{M\setminus B}\big)
  \quad\mbox{ in $L^p(M\setminus
    B,\wedge^{2j}T^*M)\ \forall 1<p<\tfrac nj$   }
\end{equation*}
as $i\to\infty$,
provided $1\le j\le n-1$ and $n>2$. In view of
Lemma~\ref{lem:closed_cohom_class}, this convergence
and~(\ref{Chern-i}) imply the claim~(\ref{Chern-claim}).

Finally, in the case of the Yang-Mills equation $n=2$, the weak
convergence $F_{A_i^\alpha}\wto F_{A^\alpha}$ in $L^2$
and the linearity of $p_1$ imply 
\begin{equation*}
  p_1\big(\tfrac{\mathbf{i}}{2\pi}F_{A_i}\big|_{M\setminus B}\big)
  \wto
  p_1\big(\tfrac{\mathbf{i}}{2\pi}F_{A}\big|_{M\setminus B}\big)
  \quad\mbox{ weakly in $L^2(M\setminus
    B,\wedge^2T^*M)$   }
\end{equation*}
as $i\to\infty$. This implies by Lemma~\ref{lem:closed_cohom_class} 
that the first Chern class is preserved in the limit, which is the
only assertion claimed in~(\ref{Chern-claim}) for $n=2$. We have
thereby established the claim~(\ref{Claim}) and 
conclude that the limit bundle $\widehat P$ possesses the desired Chern
classes, which means $\widehat P\in\mathcal{P}_*(M,P_0)$. 
\\
 
\textbf{Final step: Minimization property of the limit connection.}
We choose a partition of unity $\{\zeta_\alpha\}_{\alpha\in\N}$ 
subordinate to the cover $U_\alpha$ of $M\setminus\Sigma$. 
From~(\ref{F-converge}) and the lower semi-continuity of the norm with
respect to weak convergence, we infer
\begin{align*}
  \int\zeta_\alpha \big(|\das{n-2}F_A|^2+|F_A|^n\big)\dx
   &\le
  \lim_{i\to\infty}\int\zeta_\alpha
  \big(|\dasi{n-2}F_{A_i}|^2+|F_{A_i}|^n\big)\dx\\
  &=\lim_{i\to\infty}\int\zeta_\alpha\,d\mu_i=\int\zeta_\alpha\,d\mu,
\end{align*}
where we used the weak*-convergence $\mu_i\wsto\mu$ in the last step.
Summing over $\alpha\in\N$ and using the fact that $\zeta_\alpha$ is a
partition of unity on $M\setminus\Sigma$, we arrive at 
\begin{align*}
  Y_n(\widehat A,\widehat P)
  &=
  Y_n(A,P)
  \le 
  \mu(M\setminus\Sigma)
  \le
  \mu(M)\\
  &=
  \lim_{i\to\infty}\mu_i(M)
  =
  \lim_{i\to\infty}Y_n(A_i,P_i)
  =\inf Y_n.
\end{align*}
This implies that the pair $(\widehat A,\widehat P)$ minimizes the functional 
$Y_n$ in the class of principal bundles $P\in\mathcal{P}_*(M,P_0)$ and of
smooth connections $A$ on $P$. 
\qed  

\subsection{The subcritical dimensions}
\label{sec:existence-subcritical}
If $m<2n$, we can even minimize $Y_n$ in an arbitrary fixed 
bundle over $M$:

\begin{theorem}\label{thm:low_existence}
  Let $M$ be a compact manifold of dimension $m<2n$ and $P$ a smooth
  principal $G$-bundle over $M$ with compact structure group $G$.
  Then there is a smooth connection $A$ minimizing $Y_n$ (or $Z_n$) among
  all connections on $P$ of class $W^{n-1,2}$.
\end{theorem}

In order to show that the bundle does not change when passing to the
limit in a minimizing sequence, we need the following patching
construction that goes back to Uhlenbeck \cite[Prop. 3.2]{Uh2}. A detailed proof
can be found in \cite[Lemma 7.2(i)]{We}, where the result is even
proven for transition
functions $g_{\alpha\beta}, h_{\alpha\beta}$ that are only of class $W^{k+1,p}$ for
any $k\in\N$ and $p>\frac m2$, with resulting gauge transformations of
the same class. Since we can choose $k\in\N$ arbitrarily, this result
includes the case of $C^\infty$ functions.

\begin{lemma}\label{lem:patching}
  Consider a locally
  finite open cover $M=\bigcup_{\alpha\in\N}U_\alpha$ of a compact manifold
  $M$ by precompact sets $U_\alpha$. There are a constant $\delta_M$,
  only depending on the geometry of $M$, and open subsets
  $V_\alpha\subset U_\alpha$, $\alpha\in\N$,
  still covering $M$ such that the following
  holds. For any two sets of transition functions
  $g_{\alpha\beta},h_{\alpha\beta}\in C^\infty(U_\alpha\cap U_\beta,
  G)$, $\alpha,\beta\in\N$, that both satisfy the cocycle conditions 
  \begin{equation*}
    g_{\alpha\beta}g_{\beta\gamma}=g_{\alpha\gamma}
    \quad\mbox{and}\quad
    h_{\alpha\beta}h_{\beta\gamma}=h_{\alpha\gamma}
    \qquad\mbox{on }U_\alpha\cap U_\beta\cap U_\gamma
  \end{equation*}
   and furthermore
  \begin{equation}\label{Linfty-close}
    \sup_{x\in U_\alpha\cap U_\beta}
    \mathrm{d}(g_{\alpha\beta}(x),h_{\alpha\beta}(x))\le \delta_M
    \qquad\mbox{for all } \alpha,\beta\in\N,
  \end{equation}
  there exist local gauge transformations $h_\alpha\in
  C^\infty(V_\alpha,G)$ for all $\alpha\in\N$ such that
  \begin{equation*}
    h_\alpha^{-1}h_{\alpha\beta}h_\beta=g_{\alpha\beta}
   \qquad\mbox{on }V_\alpha\cap V_\beta,
  \end{equation*}
  provided $V_\alpha\cap V_\beta\neq\varnothing$.
\end{lemma}

After these preparations, we can proceed to the \\

{\bf Proof of Theorem \ref{thm:low_existence}.} 
We start by choosing a minimizing sequence of $W^{n-1,2}$ connections
$A_i$, $i\in\N$ on the bundle $P$. Theorem~\ref{WA-estimates} provides
us with the bound
\begin{equation}\label{sup-WA-Bound}
  \sup_{i\in\N}\,\|D_{A_i}^\ell F_{A_i}\|_{L^{2n/(\ell+2)}(M)}<\infty
  \qquad\mbox{for }\ell=0,\ldots,n-2.
\end{equation}

\textbf{Step 1: Construction of local Uhlenbeck gauges.}
Since we are in the subcritical dimension $m<2n$, we can
always choose a radius $r_0>0$ so small that
\begin{equation}
  \label{no-concentration}
  \sup_{i\in\N}\ 
  r^{2n-m}\int_{B_r}
  |F_{A_i}|^n\,dx
  \le 
  r^{2n-m}\sup_{i\in\N}Y_n(A_i)< \kappa^n
\end{equation}
holds for any ball $B_r\subset M$ of radius $r\le r_0$, where here,
we chose the constant $\kappa=\kappa(M)>0$ from the Gauge
Theorem~\ref{uhl_n}. 
Now we cover $M$ by finitely many balls $U_\alpha$,
$\alpha=1,\ldots,L$ of the same radius $r\le \frac12 r_0$ 
in such a way that the bundle
$P$ is trivial over $U_\alpha$ for each $\alpha=1,\ldots,L$. 
This means that we can find bundle isomorphisms
\begin{equation*}
  U_\alpha\times G\to P|_{U_\alpha},\quad 
  (x,g)\mapsto g_\alpha(x) g
\end{equation*}
given by $C^\infty$-sections $g_\alpha:U_\alpha\to P$. 
The local trivializations give rise
to transition functions $g_{\alpha\beta}\in C^\infty(U_\alpha\cap
U_\beta,G)$ defined by
\begin{equation}\label{bundle-transition}
  g_\beta=g_\alpha\,g_{\alpha\beta}\qquad\mbox{on } U_\alpha\cap U_\beta
\end{equation}
for all $\alpha,\beta$ for which the intersection $U_\alpha\cap
U_\beta$ is non-empty. From~(\ref{no-concentration}), we know that 
\begin{equation}
  \label{no-concentration-U}
  \sup_{i\in\N}\
  (2r)^{2n-m}\int_{2U_\alpha}
  |F_{A_i}|^n\,dx
  <
  \kappa^n.
\end{equation}
The Gauge Theorem~\ref{uhl_n} thus 
yields gauge transformations $u_{i,\alpha}\in W^{n,2}(U_\alpha,G)$ such that the
localized connection forms 
$A_i^{\alpha}:=u_{i,\alpha}^*A_i$ satisfy
$d^*A_i^{\alpha}=0$ and
\begin{align}\label{A-Bounds}
  &\sup_{i\in\N} \sum_{\ell=0}^{n-1}r^{\frac{2n-m}{2n}(\ell+1)}\|D^\ell
  A_i^{\alpha}\|_{L^{2n/(\ell+1)}(U_\alpha)}\\\nonumber
  &\qquad\le
  C\sup_{i\in\N}
   \Big(r^{\frac{2n-m}{2}}\|D_{A_i}^{n-2}F_{A_i}\|_{L^2(2U_\alpha)}
    +r^{\frac{2n-m}{n}}\|F_{A_i}\|_{L^n(2U_\alpha)}\Big)
  <\infty
\end{align}
for every $\alpha=1,\ldots,L$, where the finiteness of the right-hand
side follows from
(\ref{sup-WA-Bound}). The gauge transformations $u_{i,\alpha}$ define transition functions
$\phi^i_{\alpha\beta}\in W^{n,2}(U_\alpha\cap U_\beta,G)$ by the identity
\begin{equation*}
  u_{i,\beta}=u_{i,\alpha}\phi^i_{\alpha\beta}
  \qquad\mbox{on }U_\alpha\cap U_\beta.
\end{equation*}
From $A_i^\beta=(\phi^i_{\alpha\beta})^*A_i^\alpha$ we infer
\begin{equation*}
  d\phi^i_{\alpha\beta}
  =
  \phi_{\alpha\beta}^iA_i^\beta-A_i^\alpha\phi_{\alpha\beta}^i
  \qquad\mbox{on }U_\alpha\cap U_\beta.
\end{equation*}

\textbf{Step 2: Convergence to a smooth limit connection. }
This step is analogous to the critical case in
Theorem~\ref{thm:existence}. 
Starting from the estimate~(\ref{A-Bounds}), we can argue as in the
critical case in order to find local
limit connections $A_*^\alpha\in
W^{n-1,2}(U_\alpha,\wedge^1\R^m\otimes\gg)$ and transition functions
$\phi_{\alpha\beta}\in W^{n,2}(U_\alpha\cap U_\beta, G)$ such that
after extraction of a subsequence we
have  
\begin{equation}\label{A-phi-converge}
  \left\{
  \begin{array}{c l}
  A_i^{\alpha}\wto A_*^{\alpha}&\mbox{weakly in
  }W^{n-1,2}(U_\alpha,\wedge^1\R^m\otimes\gg),
  \\[1ex]
 \phi^i_{\alpha\beta}\wto \phi_{\alpha\beta}
 & \mbox{weakly in }W^{n,2}(U_\alpha\cap U_\beta,G).
\end{array}
\right.
\end{equation}
The weak convergence implies that the local limit connections are
still in Uhlenbeck gauge, i.\,e. 
\begin{equation}\label{Uhlenbeck-limit-2}
  d^*A_*^{\alpha}=0\qquad\mbox{on }U_\alpha.
\end{equation}
Moreover, the limit connections satisfy the transformation rule
\begin{equation}
  \label{transition-limit-2}
  d\phi_{\alpha\beta}
  =
  \phi_{\alpha\beta}A_*^{\beta}-A_*^{\alpha}\phi_{\alpha\beta}
  \qquad\mbox{on }U_\alpha\cap U_\beta,
\end{equation}
and the cocycle conditions are also preserved in the limit in the sense
\begin{equation}\label{cocycle-2}
  \phi_{\alpha\beta}\phi_{\beta\gamma}=\phi_{\alpha\gamma}
  \qquad\mbox{on }U_\alpha\cap U_\beta\cap U_\gamma.
\end{equation}
As in step 2 of the proof of Theorem~\ref{thm:existence} we deduce
that $A_*^\alpha$ satisfies the Euler-Lagrange system for $Y_n$. The
proof of this fact becomes only easier in the subcritical case because now
the compact embedding $W^{n-1,2}\hookrightarrow W^{\ell,2n/\ell+1}$
for $\ell<n-1$ implies
strong convergence of the lower order terms $D^\ell {A_i^\alpha}\to
D^\ell {A_*^\alpha}$ in $L^{2n/\ell+1}$. Having established the Euler-Lagrange system and the Coulomb
condition~(\ref{Uhlenbeck-limit-2}) we can apply the Regularity
Theorem~\ref{thm:regularity} in order to deduce $A_*^\alpha\in
C^\infty(U_\alpha)$ for every $\alpha=1,\ldots,L$. From the
transformation rule~(\ref{transition-limit-2}) we then also infer
$\phi_{\alpha\beta}\in C^\infty(U_\alpha\cap U_\beta)$ for 
$\alpha,\beta=1,\ldots,L$.   
From~(\ref{transition-limit-2}) and~(\ref{cocycle-2}) we know that the
local connection forms $A_*^\alpha$ stem from a smooth connection
$A_*$ on a new smooth principal bundle $P_*$ that is given
by the transition functions $\phi_{\alpha\beta}$. \\

\textbf{Step 3: Identification of the limit bundle and conclusion.}
It remains to show that the bundle $P_*$, determined by the transition functions
$\phi_{\alpha\beta}$, is isomorphic to the original bundle $P$ with
transition functions $g_{\alpha\beta}$. This step relies crucially on
the subcritical dimension $m<2n$, in particular on the compact embedding
$W^{n,2}\hookrightarrow C^0$ that fails to hold in the critical
dimension. 
In view of this embedding and the weak
convergence~(\ref{A-phi-converge})$_2$, we can assume that all
$\phi^i_{\alpha\beta}$ are uniformly close to $\phi_{\alpha\beta}$ in
the sense
\begin{equation}
  \label{phi-uniform-conv}
  \sup_{x\in U_\alpha\cap
    U_\beta}\mathrm{d}\big(\phi_{\alpha\beta}^i(x),\phi_{\alpha\beta}(x)\big)
  \le
  \tfrac12 \delta_M
  \qquad\mbox{for all }i\in\N, 
\end{equation}
where $\delta_M$ denotes the constant from Lemma~\ref{lem:patching}. 
Next, we observe that
\begin{equation*}
  \phi^1_{\alpha\beta}=\big(g_\alpha^{-1}u_{1,\alpha}\big)^{-1}
  g_{\alpha\beta} \big(g_\beta^{-1}u_{1,\beta}\big)
  \qquad\mbox{on }U_\alpha\cap U_\beta. 
\end{equation*}
The gauge transformations $g_\alpha^{-1}u_{1,\alpha}$ may not be smooth,
but by the embedding $W^{n,2}\hookrightarrow C^0$ they are
continuous. Therefore they can be approximated uniformly by smooth
gauge transformations $v_\alpha\in C^\infty(U_\alpha)$. Letting
$\tilde g_{\alpha\beta}:= v_\alpha^{-1}g_{\alpha\beta}v_\beta$, we can
thereby achieve
\begin{equation*}
   \sup_{x\in U_\alpha\cap
    U_\beta}\mathrm{d}\big(\tilde g_{\alpha\beta}(x),\phi^1_{\alpha\beta}(x)\big)
  \le
  \tfrac12 \delta_M.
\end{equation*}
Joining this with~(\ref{phi-uniform-conv}), we conclude
\begin{equation*}
  \sup_{x\in U_\alpha\cap
    U_\beta}\mathrm{d}\big(\tilde g_{\alpha\beta}(x),\phi_{\alpha\beta}(x)\big)
  \le
  \delta_M.
\end{equation*}
It is straightforward to check that
the new transition functions $\tilde g_{\alpha\beta}$ still
satisfy the cocycle conditions. We are therefore in a position to
apply Lemma~\ref{lem:patching}, which provides us with open sets
$V_\alpha\subset U_\alpha$ that still cover $M$ and gauge
transformations $\phi_\alpha\in C^\infty(V_\alpha,G)$ with
\begin{equation*}
  \phi_\alpha^{-1}\phi_{\alpha\beta}\phi_\beta
  =
  \tilde g_{\alpha\beta}
  =
  v_\alpha^{-1}g_{\alpha\beta}v_\beta
  \qquad\mbox{on }V_\alpha\cap V_\beta. 
\end{equation*}
Abbreviating $\psi_\alpha:=\phi_\alpha v_\alpha^{-1}$, we can re-write
this to 
\begin{equation*}
  \psi_\alpha^{-1}\phi_{\alpha\beta}\psi_\beta
  =
  g_{\alpha\beta}
  \qquad\mbox{on }V_\alpha\cap V_\beta. 
\end{equation*}
This means that $P_*$ is isomorphic to $P$
by the bundle isomorphism that is locally given by $\psi_\alpha:=
\phi_\alpha v_\alpha^{-1}\in C^\infty(V_\alpha,G)$. 
Letting $A^\alpha:=\psi_\alpha^*A_*^\alpha$, we observe
\begin{equation*}
  g_{\alpha\beta}^*A^\alpha = (\psi_\alpha g_{\alpha\beta})^*
  A_*^\alpha
  =(\phi_{\alpha\beta}\psi_\beta)^*A_*^\alpha=A^\beta.
\end{equation*}
Therefore, the $A^\alpha$ are the localizations of a connection $A$
on the bundle $P$. In order to check that this is the desired
minimizer, we use the lower semicontinuity of the norm with respect
to weak convergence for the estimate
\begin{align*}
  Y_n(A)=Y_n(A_*)\le \lim_{i\to\infty}Y_n(A_i)=\inf Y_n,
\end{align*}
where the infimum is taken over all $W^{n-1,2}$-connections on
$P$. This shows the minimization property of $A$ and concludes the
proof of the theorem.
\qed

\vfill

\begin{center}
\scriptsize
Fakult\"at f\"ur Mathematik, Universit\"at Duisburg-Essen, D-45117 Essen, 
Germany.\\
{\tt andreas.gastel@uni-due.de}\\
{\tt christoph.scheven@uni-due.de}
\end{center}

\end{document}